\documentclass{asl}

\newcommand{\smn}{S_{m,n}}
\newcommand{\Smn}{S_{m+M,n+N}}
\newcommand{\Dom}{K\text{-Dom}_{m+M,n+N}}
\newcommand{\dom}{K\text{-Dom}_{m,n}}
\newcommand{\kDom}{K'\text{-Dom}_{m+M,n+N}}
\newcommand{\kdom}{K'\text{-Dom}_{m,n}}
\newcommand{\ls}{\left<}
\newcommand{\rs}{\right>}
\newcommand{\lb}{\left[\left[ }
\newcommand{\rb}{\right]\right]_s}
\newcommand{\pb}{\bar{p}}
\newcommand{\qb}{\bar{q}}
\newcommand{\mfp}{\mathfrak{p}}
\newcommand{\mfq}{\mathfrak{q}}
\newcommand{\mfm}{\mathfrak{m}}
\newcommand{\mfn}{\mathfrak{n}}

\theoremstyle{plain}
\newtheorem{thm}{Theorem}[section]
\newtheorem{lem}[thm]{Lemma}
\newtheorem{cor}[thm]{Corollary}

\theoremstyle{definition} 
\newtheorem{defn}[thm]{Definition}

\newtheorem{rem}[thm]{Remark}

\numberwithin{equation}{section}

\begin{document}

\author{Y. F\i rat \c{C}el\.{i}kler}
\revauthor{\c{C}el\.ikler, Y. F\i rat}

\title[Dimension theory and normalization for $D$-semianalytic sets]{Dimension theory and parameterized normalization for $D$-semianalytic sets over non-Archimedean fields}
\address{Mathematics Department\\ Purdue University\\ West Lafayette\\ IN 47907, USA.}

\email{ycelikle@math.purdue.edu}

\thanks{This paper is a part of the author's  PhD thesis and the author wishes to thank his advisor Leonard Lipshitz for support and guidance during the preparation of this paper. He also wishes to thank William Heinzer and Zachary Robinson for  helpful conversations and the referee for pointing out significant improvements in an earlier version of this paper.}
\thanks{The author was supported in part by NSF grant number DMS 0070724.}

\subjclass[2000]{Primary (32B20); Secondary (03C10), (12J25)}

\begin{abstract}
We develop a dimension theory for $D$-semianalytic sets over an arbitrary non-Archimedean complete field. Our main results are the equivalence of several notions of dimension and a theorem on additivity of dimensions of projections and fibers in characteristic $0$. We also prove a parameterized version of normalization for $D$-semianalytic sets.
\end{abstract}

\maketitle

\section{Introduction and Notation}

In this paper, we develop and investigate several algebraic and geometric notions of dimension (see Definitions \ref{Adim}, \ref{gdim} and \ref{wdim}) which are associated with $D$-semianalytic sets over an arbitrary non-Archimedean complete field $K$. Our main result (Theorem \ref{alleq0}) is the equivalence of these notions of dimension in characteristic $0$. We also prove some weaker relations between them in characteristic $p>0$. In addition, a theorem on additivity of dimensions of fibers and projection (Theorem \ref{dimadd}) is proved in characteristic $0$. One of the main tools for proving these results is the Parameterized Normalization Lemma (Lemma \ref{norm1}) which also has applications to other problems concerning the geometry of $D$-semianalytic sets apart from the ones mentioned in this paper. The $D$-semianalytic sets that we are considering in this paper are the subsets of $(K^{\circ})^m\times (K^{\circ\circ})^n$  (see Definition \ref{k0}) which are quantifier free definable in the three sorted language which takes the members of $S_{m,n}(E,K)$ (see Definition \ref{separ}) as the analytic functions over $(K^{\circ})^m\times (K^{\circ\circ})^n$ as defined by Lipshitz in \cite{rigid}. Nevertheless we will avoid discussing this language in detail and work with the Definition \ref{Dsemi} for $D$-semianalytic sets instead.

 In \cite{strat} Lipshitz and Robinson proved a smooth stratification theorem and obtained some results similar to the ones in this paper in dimension theory of quasi-affinoid subanalytic subsets of $(\bar{K}^{\circ})^m\times (\bar{K}^{\circ\circ})^n$ where $K$ is an arbitrary non-Archimedean complete field and $\bar{K}$ is an algebraically closed, complete extension of $K$. These are definable (with quantifiers) subsets in the same language as above containing members of $S_{m,n}(E,K)$ as the analytic functions. Let $K'$ be a complete extension of a complete non-Archimedean field $K$, let $\mathcal{L}_1$ and $\mathcal{L}_2$ be two such languages having members of $\smn(E,K)$ and $\smn(E,K')$ as their analytic functions respectively. As a subset of $({K'}^{\circ})^m\times({K'}^{\circ\circ})^n$ which is quantifier free definable in $\mathcal{L}_1$ is quantifier free definable in $\mathcal{L}_2$, the distinction between the coefficient field and the field in which the points lie is unnecessary in our case.

 We also would like to mention some related results of other authors. Denef and van den Dries introduced $p$-adic subanalytic sets and developed a dimension theory for them in \cite{denefvdd}.  By the Quantifier Elimination Theorem (Theorem 1.1) of \cite{denefvdd} these subanalytic sets are in fact quantifier free definable subsets of $\mathbb{Q}_p^m$ in the language containing function symbols for each member of the Tate algebra $T_m$ for all $m$, a function symbol $D$ for restricted division and  the $n^{\text{th}}$-power predicates $P_n$ as defined in \cite{macint} by Macintyre. Therefore they are analogues of $D$-semianalytic sets in our setting. In \cite{pmin} Haskell and Macpherson developed a dimension theory for definable sets in the $p$-minimal setting. On the other hand, in \cite{dimdef} van den Dries considers dimension functions over Tarski systems of definable sets and proves existence of a unique dimension function (algebraic dimension) for Tarski system of definable sets of each algebraically bounded henselian field of characteristic $0$. There are several similarities between his and our results in spite of the differences between the settings. The main difference is that we are working on {\em quantifier free definable sets} in a specific language containing analytic functions whereas in \cite{dimdef}, the author is interested in {\em all definable sets} in any suitable language provided that the structure is algebraically bounded. Another difference is that the key concepts of \cite{dimdef} like algebraic dimension and algebraic boundedness, as they are stated, behave well for sets definable in an {\em algebraic language} as opposed to an {\em analytic language} like ours. Nevertheless, the Parameterized Normalization Lemma gives a weaker version of ``algebraic boundedness'' which is good enough to carry out the somewhat similar arguments in our analytic setting. In a similar way the restricted Krull dimension (Definition \ref{Adim}) can be seen as an analytic analogue of algebraic dimension of \cite{dimdef} and our main results Theorem \ref{alleq0} and Theorem \ref{dimadd} state that, in case {\em Char} $K=0$, the notions of dimensions that we consider in this paper are all equivalent and satisfy conditions similar to the axioms of a dimension function of \cite{dimdef}. As those axioms are stated for a dimension function over a Tarski system of sets, some minor adjustments are needed. 

Our main motivation and basic tools come from the following sources. In \cite{rigid} it was shown that the theory of algebraically closed complete non-Archimedean fields in the three sorted language involving the members of the rings of separated power series $\smn$, the norm function $|\cdot|$  and the restricted division operations $D_0$ and $D_1$ admits quantifier elimination. These division operations are as follows.
$$\begin{array}{rcl}
D_0:K^2\rightarrow K^{\circ},\text{ defined as } D_0(x,y)&=&\begin{cases}x/y &\text{if } |x|\leq|y|\neq 0\\0&\text{otherwise}\end{cases}\\
D_1:K^2\rightarrow K^{\circ\circ},\text{ defined as } D_1(x,y)&=&\begin{cases}x/y& \text{if } |x|<|y|\\0&\text{otherwise.}\end{cases}\\
\end{array}$$
As an immediate corollary to this result, one sees that the quantifier free definable and definable sets in this language over algebraically closed complete non-Archimedean fields coincide. Lipshitz and Robinson generalized this result to more general classes of definable sets in \cite{modcomp}. In \cite{strat} some of the geometric properties of the subanalytic subsets of $(\bar{K}^{\circ})^m\times (\bar{K}^{\circ\circ})^n$  were investigated where $\bar{K}$ is an algebraically closed, complete non-Archimedean field and we will use some of the  same ideas. On the other hand Lipshitz and Robinson established many results in commutative algebra of the rings of separated power series and related algebras in \cite{RoSPS}. Results from this paper will help us obtain algebraic geometry like results between the $D$-semianalytic sets and the algebras associated with them. Noetherianness of (a property that we  use without mentioning throughout the paper), and having Weierstrass Division Theorems for the separated power series rings are two of the key results we will often borrow from this source.

In more detail, the outline of this paper is as follows. After preliminary ground work, in Section 3, we introduce a restricted Krull dimension of algebras associated with $D$-semianalytic sets and prove its nice behavior under operations like adjoining fractions and ground field extensions. In Section 4 we investigate a geometric notion of dimension and prove that, in case we are working over a field of characteristic $0$, a $D$-semianalytic set can be decomposed into finitely many $D$-semianalytic manifolds. The maximum of geometric dimensions of such manifolds will dominate the restricted dimension of the associated algebra. Next we turn our attention to a weaker notion of geometric dimension and make a first attempt to establish the relations between these dimensions. In Section 5 it is proved (Lemma \ref{norm1}) that given a $D$-semianalytic set $X$, it is possible to obtain a finite collection of normalized quasi-affinoid algebras the union of whose associated $D$-semianalytic sets is $X$. Furthermore if some of the variables in the algebra are designated to be parameter variables then the normalization maps preserve the parameter structure in the algebra and the restricted dimension as well. As in many of the results mentioned above, normalization is fundamental in establishing the link between the geometry of definable sets and the algebras of functions over them. As another application, apart from the ones mentioned in this paper, one can easily obtain the Quantifier Elimination Theorem (Theorem 3.8.1) of \cite{rigid} using the Parameterized Normalization Lemma. Finally in Section 6 we make use of this normalization to establish further connections between the notions of dimension. Lastly we prove (Theorem \ref{dimadd}) that, in characteristic $0$, if a $D$-semianalytic set's projection onto a subspace of dimension $m$ contains a dense subset of points with fibers of dimension $d$, then the dimension of the $D$-semianalytic set is at least $d+m$.

Throughout this paper $K$ will denote a field which is complete with respect to a non-Archimedean norm, $K'$ will denote an arbitrary complete field extension of $K$ and $\bar{K}$ will denote an algebraically closed complete extension of $K$.

\begin{defn}\label{k0}
For any non-Archimedean complete field $K$, $K^{\circ}$ will denote its valuation ring and $K^{\circ\circ}$ will denote the maximal ideal of $K^\circ$. Or in other words
$$\begin{array}{rcl}
K^{\circ}&=&\{a\in K:|a|\leq 1\},\\
K^{\circ\circ}&=&\{a\in K:|a|< 1\}.
\end{array}$$
\end{defn}

We are interested in analytic functions over the set $(K^{\circ})^m\times(K^{\circ\circ})^n$ and a natural choice for such functions is the members of separated power series rings $\smn(E,K)$. Following Definition 2.1.1 of \cite{RoSPS}, these rings are obtained as follows.

\begin{defn}\label{separ}
Let $x=(x_1,...,x_m)$ and $\rho=(\rho_1,...,\rho_n)$ denote multi variables, fix a complete, quasi-Noetherian subring $E$ of $K^{\circ}$ (which also has to be a Discrete valuation ring in case $Char\: K=p>0$) and let $\{a_i\}_{i\in\mathbb{N}}$ be a zero sequence in $K^\circ$, and $B$ be the  local quasi-Noetherian ring
$$(E[a_0,a_1,...]_{\{a\in E[a_0,a_1,...]:|a|=1\}})^{\wedge},$$
where $\text{}^\wedge$ denotes the completion in $|\cdot|$. Let $\mathcal{B}$ be the family of all such rings. Define the separated power series ring
$$\smn(E,K):=K\otimes_{K^\circ}\left(\lim_{\overrightarrow{B\in\mathcal{B}}}B\ls x\rs[[\rho]]\right).$$

For $f=\sum_{\alpha,\beta}a_{\alpha,\beta}x^{\alpha}\rho^{\beta}\in\smn(E,K)$, the {\em Gauss norm} of $f$ is defined as 
$$||f||:=\sup_{\alpha,\beta}|a_{\alpha,\beta}|.$$
\end{defn} 

Note that in general the rings $\smn(E,K)$ are not complete in the Gauss norm, but for many choices of the quasi-affinoid ring $E$ they are complete. We refer the reader to Remark 2.1.2 and Theorem 2.1.3 of \cite{RoSPS} for a detailed discussion of this aspect of the $\smn(E,K)$.

We will usually denote $\smn(E,K)$ by $\smn$, which will not lead to a confusion. The elements of a separated power series ring are convergent on the set $(K^{\circ})^m\times(K^{\circ\circ})^n$. The term {\em quasi-affinoid} is used to refer to objects associated with such rings. 

In order to incorporate the restricted division operations $D_0$ and $D_1$ into an algebraic context, we follow \cite{strat} and work in generalized rings of fractions. To be able to define those rings, we first introduce more terminology. A {\em quasi-affinoid algebra} $A$ is an algebra of the form $A:=\smn/J$ for some integers $m$, $n$ and ideal $J$ of $\smn$. For a quasi-affinoid algebra $A=\smn/J$, variables $y=(y_1,...,y_M)$ and $\lambda=(\lambda_1,...,\lambda_N)$ not appearing in $\smn$, define
$$A\ls y\rs\lb \lambda\rb:=\Smn/J\Smn.$$
More generally if $A$, $B_1$ and  $B_2$ are quasi-affinoid algebras, where 
$$B_1=A\ls y_1,....,y_{m_1}\rs\lb\lambda_1,...,\lambda_{n_1}\rb/I_1$$ 
and
$$B_2=A\ls y_{m_1+1},....,y_{m_2}\rs\lb\lambda_{n_1+1},...,\lambda_{n_2}\rb/I_2$$
then as in Definition 5.4.2 of \cite{RoSPS}, define the separated tensor product of $B_1$ and $B_2$ over $A$ as
$$B_1\otimes_A^sB_2:=A\ls y_1,...,y_{m_2}\rs\lb \lambda_1,...,\lambda_{n_2}\rb/(I_1\cup I_2).$$
By Theorem 5.2.6 of \cite{RoSPS} this product is independent of presentations of $B_1$ and $B_2$. 
 
Let $K'$ be a complete field extension of $K$ and $E'\subset K'$ be a complete, quasi-Noetherian ring. Assume $\smn(E,K)\subset\smn(E',K')$ and let $B_1$ be as above, then we will follow Definition 5.4.9 of \cite{RoSPS} and say the $K'$-affinoid algebra
$$B_1'=S_{0,0}\otimes^s_{S_{0,0}(E,K)}B_1:=\smn(E',K')/I_1\cdot\smn(E',K')$$
results from $B_1$ by ground field extension from $(E,K)$ to $(E',K')$. Note that by Lemma 4.2.8 of \cite{RoSPS}, $\smn(E',K')$ is faithfully flat over $\smn(E,K)$.

With this notation established, we are ready to define one of the main objects of our study.

\begin{defn}\label{grf}
Following the Definition 5.3.1 of \cite{RoSPS}, a {\em generalized ring of fractions} over a quasi-affinoid algebra $A$ is inductively defined as follows:

i) $A$ is a generalized ring of fractions over $A$.

ii) If $B$ is a generalized ring of fractions over $A$, $y$, $\lambda$ multi-variables not appearing in the presentation of $B$, and $g,\:f_1,...,f_M,F_1,...,F_N\in B$ , then 
$$B\ls f/g\rs\lb F/g\rb:=B\ls y\rs\lb\lambda\rb/(\{gy_i-f_i\}_{i=1}^M\cup\{g\lambda_j-F_j\}_{j=1}^N)$$
is also a generalized ring of fractions over $A$.
\end{defn}

We would like to note that this definition is a generalization of the definition of affinoid generalized rings of fractions. These rings are discussed in Subsection 6.1.4 of \cite{bgr}.

Let
\begin{equation}\label{bova}\tag{$\ast$}
B=A\ls y_1,...,y_M\rs\lb \lambda_1,...,\lambda_N\rb/(\{g_iy_i-f_i\}_{i=1}^M\cup\{G_j\lambda_j-F_j\}_{j=1}^N)
\end{equation}
be a generalized ring of fractions over the quasi affinoid algebra $A=\smn/J$ where $g_i$, $f_i$, $G_j$, $F_j$ are elements of $\Smn$. There is a natural correspondence between a maximal ideal $\mfm$ of $\Smn$ which contains the ideal $J\cup (\{g_iy_i-f_i\}_{i=1}^M\cup\{G_j\lambda_j-F_j\}_{j=1}^N)$ and a maximal ideal $\mfm_B$ of $B$. With this notation we can state the next definition.

\begin{defn}\label{dom}
Let $B$ be as in Equation \eqref{bova}, then the {\em domain of $B$ over $A$ } is the set
\begin{multline*}
\text{Dom}_AB:=\{\mfm_B\in \text{Max }B:J\cup\{g_iy_i-f_i\}_{i=1}^M\cup\{G_j\lambda_j-F_j\}_{j=1}^N\subset \mfm\\
\text{ and }g_i,\:G_j\not\in\mfm\text{ for all }i,j\}.
\end{multline*}
\end{defn}

\begin{rem}\label{domr}

i) This definition depends on the specific presentation of the quasi-affinoid algebra $B$, but is independent of the representatives chosen for $g_i$, $f_i$, $G_j$, and $F_j$ in the presentation.

ii) The reader should be warned that above definition is slightly different than the definition of the domain in Definition 2.2 of \cite{modcomp} which is a subset of Max $A$ rather than Max $B$. We choose to make this change for simplicity in later arguments. Actually by Theorem 4.1.1 (Nullstellensatz) of \cite{RoSPS}, for each maximal ideal $\mfm$ of $\Smn$, $\Smn/\mfm$ is an algebraic extension of $K$ and by this theorem and by induction it is easy to see that there is a natural one to one (but not necessarily onto) map $\sigma_B:\text{ Dom}_AB\rightarrow  \text{ Max }A$. The image of the map $\sigma_B$ is the definition for Dom$_AB$ in \cite{modcomp}. Therefore the deviation from the notation of other authors is not significant.

iii) Let $A$ be a quasi-affinoid algebra and $B_1$ and $B_2$ be two isomorphic quasi-affinoid $A$ algebras which are generalized rings of fractions over $A$, and let $\sigma_{B_1}:\text{Dom}_{A}B_1\rightarrow \text{Max }A$ and $\sigma_{B_2}:\text{Dom}_{A}B_2\rightarrow \text{Max }A$ be the maps mentioned above, then $\sigma_{B_1}(\text{Dom}_{A}B_1)$ and $\sigma_{B_2}(\text{Dom}_{A}B_2)$ coincide by Remark 2.3 of \cite{modcomp}. 
\end{rem}

Let $B$ be a generalized ring of fractions over $\smn$. Since Dom$_AB$ is defined in terms of maximal ideals, it is more closely related to quantifier free definable subsets of $(\bar{K}^{\circ})^m\times(\bar{K}^{\circ\circ})^n$ than those of $(K^\circ)^{m}\times(K^{\circ\circ})^{n}$ where $\bar{K}$ is an algebraically closed complete extension of $K$. As we are interested in quantifier free definable subsets of $(K^\circ)^{m}\times(K^{\circ\circ})^{n}$ we introduce another notation to denote the $K'$-rational points in domains over $\smn$ where $K'$ is a complete extension of $K$. 

\begin{defn}
Let $B$ be as in Equation \eqref{bova} with $A=\smn$, the {\em $K'$-rational points in Dom$_{\smn}B$}  is the projection of the set
\begin{multline*}
\{\pb \in ({K'}^\circ)^{m+M}\times({K'}^{\circ\circ})^{n+N}: (g_iy_i-f_i)(\pb)=0, (G_j\lambda_j-F_j)(\pb)=0\\ 
\text{ and }g_i(\pb)\neq 0, G_j(\pb)\neq 0 \text{ for all } i,j\}
\end{multline*}
onto $({K'}^\circ)^{m}\times({K'}^{\circ\circ})^{n},$ and it is denoted by $\kdom B$. 
\end{defn}

Notice that the above projection is one to one, analytic and is also the restriction of the aforementioned map $\sigma_{B_1}$ between Dom$_{\smn}B$ and Max $\smn$ to maximal ideals that correspond to points in $({K'}^\circ)^{m+M}\times({K'}^{\circ\circ})^{n+N}$. Notice also that $\kdom B$ is an open subset of $({K'}^\circ)^{m}\times({K'}^{\circ\circ})^{n}$ in the metric topology.

 If $B$ is a generalized ring of fractions over $\smn$ then in a natural way members of $B$ can be seen as analytic functions on $K'\text{-Dom}_{m,n} B$. This observation enables us to give the next definition.

\begin{defn}\label{var}
Let $B$ be a generalized ring of fractions over $\smn$ and $I$ be an ideal of $B$. We will write $K'\text{-Dom}_{m,n} B\cap V(I)_{K'}$ for the subset of $({K'}^\circ)^m\times({K'}^{\circ\circ})^n$ which consists of points in $\kdom B$ that are zeros of all the elements of $I$. 
\end{defn}

At the beginning of this paper we mentioned that the $D$-semianalytic sets are quantifier free definable subsets of $(K^\circ)^m\times(K^{\circ\circ})^n$ in the language of \cite{rigid}. This is a three sorted language containing function symbols for the members of $\smn$ for each $m$ and $n$, as well as the norm function $|\cdot|$ and the restricted division functions $D_0$ and $D_1$. By putting any quantifier free formula in this language into disjunctive normal form and observing that a negated formula of the type $\neg(f(x,\rho)=0)$ is equivalent to the formula $(|f(x,\rho)|>0)$ we arrive at the equivalent definition for $D$-semianalytic sets below.

\begin{defn}\label{Dsemi}
A {\em $D$-semianalytic subset} of $(K^\circ)^m\times(K^{\circ\circ})^n$ is a finite union of sets of the form $\dom B_i\cap V(I_i)_K$ for some generalized rings of fractions $B_i$ over $\smn$ and ideals $I_i\subset B_i$. 
\end{defn}

At this point, we would like to mention the domain of a special type of generalized ring of fractions. 

\begin{defn}\label{rdom}
If $B$ is a generalized ring of fractions over $\smn$ which is obtained by imposing the condition that at each inductive step in forming $B$, the ideal $(g,f_1,...,f_M,F_1,...,F_N)$ is the unit ideal, then Dom$_{\smn} B$ is called an {\em $R$-domain over $\smn$}.
\end{defn}

The significance of the $R$-domains comes from the fact that they are generalizations of rational domains of affinoid geometry (Definition 5 of subsection 7.2.3 of \cite{bgr}). An important property of $R$-domains is that in case $U=$ Dom$_{\smn}B$ is an $R$-domain, and $B'$ is a generalized ring of fractions over $\smn$ such that Dom$_{\smn}B'=U$, then by Proposition 5.3.6 of \cite{RoSPS} $B$ is isomorphic to $B'$. Therefore it is possible to associate a {\em quasi-affinoid ring of functions} $\mathcal{O}_K(U):=B$ to every such $R$-domain $U$. Note that although this ring of functions also depends on the choice of the quasi-Noetherian ring $E$ as well as the field $K$ (see Definition \ref{separ}), we are suppressing this dependence for simplicity of notation.

\section{Preliminaries}

 In this section we establish some basic facts about quasi-affinoid algebras which help technical aspects of our discussions in the following sections.

\begin{lem}\label{catenary}
For any maximal ideal $\mathfrak{m}$ of a quasi-affinoid algebra $A$, $A_{\mathfrak{m}}$ is a universally catenary ring.
\end{lem}

\begin{proof}
By Corollary 4.2.2 of \cite{RoSPS}, $(\smn)_{\mathfrak{m}}$ is a regular ring of Krull dimension $m+n$. By Theorem 17.8 of \cite{matsu}, any regular local ring is Cohen-Macaulay, and by Theorem 17.9 of \cite{matsu} any quotient of a Cohen-Macaulay ring is universally catenary. As $A$ is a quotient of the form $\smn/J$, the statement follows.
\end{proof}

\begin{lem}\label{flat}
$S_{m+1,n}$ and $S_{m,n+1}$ are both faithfully flat over $\smn$, hence for any quasi-affinoid algebra $A$, we have that $A\subset A\ls z\rs$ and $A\subset A\lb z\rb$ where the variable $z$ does not appear in the presentation of $A$.
\end{lem}
\begin{proof}
For the first statement, let $\bar{K}$ be an algebraically closed complete extension of $K$ and assume that we have shown that both $S_{m+1,n}(E,\bar{K})$ and $S_{m,n+1}(E,\bar{K})$ are faithfully flat over $\smn(E,\bar{K})$. By Lemma 4.2.8 of \cite{RoSPS}, $\smn(E,\bar{K})$ is faithfully flat over $\smn(E,K)$ and therefore $S_{m+1,n}(E,\bar{K})$ and $S_{m,n+1}(E,\bar{K})$ are both faithfully flat over $\smn(E,K)$ which in turn implies that $S_{m+1,n}(E,K)$ and $S_{m,n+1}(E,K)$ being faithfully flat over $\smn(E,K)$. Hence we may assume that $K$ is algebraically closed.

Next observe that by the Nullstellensatz (Theorem 4.1.1 of \cite{RoSPS}) each maximal ideal of $\smn$ is of the form $(x_1-a_1,...,x_m-a_m,\rho_1-b_1,...,\rho_n-b_n)\cdot\smn$ where $a_i,b_j\in K$, $|a_i|\leq 1$, $|b_j|<1$ for all $i,j$. Therefore by Theorem 7.2(3) of \cite{matsu}, it is enough to show that  $S_{m+1,n}$ and $S_{m,n+1}$ are flat over $\smn$ and by Theorem 7.1 of \cite{matsu}, it is enough to check that for each maximal ideal $\mfn$ of $S_{m+1,n}$ or $S_{m,n+1}$ the localizations $(S_{m+1,n})_{\mfn}$ and $(S_{m,n+1})_{\mfn}$ are flat over $(\smn)_{\mfm}$ where $\mfm=\mfn\cap\smn$. Note that again by the Nullstellensatz $\mfm$ is a maximal ideal of $\smn$.

On the other hand by Theorem 22.4 of \cite{matsu}, $(S_{m+1,n})_{\mfn}$ and $(S_{m,n+1})_{\mfn}$ are flat over $(\smn)_{\mfm}$ if and only if $((S_{m+1,n})_{\mfn})^\wedge$ and $((S_{m,n+1})_{\mfn})^\wedge$ are flat over $((\smn)_{\mfm})^\wedge$ where $^\wedge$ denotes the completion with respect to the maximal ideals. 

Now, it is easy to see that 
$$((S_{m+1,n})_{\mfn})^\wedge\simeq ((S_{m,n+1})_{\mfm})^\wedge\simeq K\left[\left[ z_1,...,z_{m+n+1}\right]\right]$$
which is faithfully flat over $K\left[\left[ z_1,...,z_{m+n}\right]\right]\simeq ((\smn)_{\mfn})^\wedge$, and the first statement of the lemma follows.

 The second statement easily follows just by writing $A=\smn/J$, $A\ls z\rs=S_{m+1,n}/JS_{m+1,n}$ (or $A\lb z\rb=S_{m,n+1}/JS_{m,n+1}$) and observing 
$$JS_{m+1,n}\cap\smn=JS_{m,n+1}\cap\smn=J$$
because of faithful flatness.
\end{proof}

As a consequence of the above lemma, for $A=\smn/J$ a quasi-affinoid algebra, $I$ an ideal of $A$ and $z$ a variable not appearing in the presentation of $A$, we get $(A/I)\ls z \rs=A\ls z \rs/IA\ls z\rs$. 

Let us make an observation about the presentation of elements in quasi-affinoid algebras before we continue with the next lemma. Suppose $\bar{f}\in A\ls z\rs$ and $f\in S_{m+1,n}$ is such that the canonical image of $f$ in $A\ls z\rs=S_{m+1,n}/JS_{m+1,n}$ is $\bar{f}$. Write $f=\sum_i c_iz^i$ for $c_i\in\smn$ and let $a_i$ denote the canonical image of $c_i$ in $A$, then $\bar{f}$ has a presentation 
$$\bar{f}=\sum_i a_i z^i.$$ 
On the other hand if $g=\sum_i d_iz^i$ is another element of $S_{m+1,n}$ whose canonical image is $\bar{f}$ and $h_1,...,h_k\in\smn$ are generators of $J$, then
$$f-g=\sum_i(c_i-d_i)z^i=h_1\sum_i a_{1i}z^i+...+h_k\sum_i a_{ki}z^i,$$
where $a_{ji}\in\smn$ for all $i$, $j$. Comparing the coefficients of $z^i$ on both sides we see that $c_i-d_i\in J$ for all $i$. This shows that such a presentation of $\bar{f}$ is unique. A similar statement is true for $A\lb z\rb$ by the same argument.

\begin{lem}\label{idealmove}
Let $A$ be a quasi-affinoid algebra, $I$ an ideal of $A$, $z$ a variable not appearing in the presentation of $A$, then 
$$IA\ls z\rs=\{\sum_{i=0}^{\infty}a_i z^i\in A\ls z\rs: a_i\in I \text{ for all } i\},$$
similarly
$$IA\lb z\rb=\{\sum_{i=0}^{\infty}a_i z^i\in A\lb z\rb: a_i\in I \text{ for all } i\}.$$
\end{lem}
\begin{proof}
For both assertions, the inclusion $\subseteq$ is clear. For the other inclusion, write $A=\smn/J$, and let $\bar{I}$ be the preimage of $I$ in $\smn$, $\bar{I}\supset J$. Let $f=\sum a_iz^i$ be an element of either $\smn\ls z\rs$ where $a_i\in \bar{I}$ for all $i$ (the case $f\in\smn\lb z\rb$ is similar). Multiplying with a constant, we may assume that $||f||=1$ where $||\cdot||$ denotes the Gauss Norm. By Lemma 3.1.7 of \cite{RoSPS}, there are $||\cdot||$-strict generators for $\bar{I}$. Let those generators be $b_1,...,b_k$. There is a $B\in\mathcal{B}$ such that for each $i$, there are $c_{1i},...,c_{ki}\in B\ls x_1,...,x_m,z\rs\lb\rho_1,...,\rho_n\rb$ which satisfy 
$$\begin{array}{rcl}
a_i&=&c_{1i}b_1+...+c_{ki}b_k\\
||a_i||&\geq& \sup_{j} \{||c_{ji}||\cdot||b_j||\}.
\end{array}$$
Hence we have a presentation of the element $f$ as
$$f=b_1\sum_i c_{1i}z^i+...+b_k\sum_i c_{ki}z^i.$$
Note that each $\sum_ic_{ji}z^i\in B\ls x_1,...,x_m,z\rs\lb\rho_1,...,\rho_n\rb$ and hence $f\in\bar{I}\smn\ls z\rs$.
\end{proof}

We wish to make a note about the contraction and extension of ideals in quasi-affinoid algebras. If $A=\smn/J$ is a quasi-affinoid algebra, it is not necessarily the case that a ring of the form $A\ls f/g\rs$ or $A\lb f/g\rb$ contains $A$, as among many other possibilities, $gz-f$ can be a unit of the ring $A\ls z\rs$ or $A\lb z \rb$. Still $A\ls f/g\rs$ and $A\lb f/g\rb$ are $A$-algebras in a natural way and so it is possible to talk about extension of an ideal $I$ of $A$ to $A\ls f/g\rs$ or $A\lb f/g\rb$ and contraction of an ideal $I'$ of $A\ls f/g\rs$ or $A\lb f/g\rb$ to $A$.

The last two lemmas in this section cover some technical aspects of the proof of the Parameterized Normalization Lemma (Lemma \ref{norm1}).

\begin{lem}\label{injective}

Let $\phi:A\hookrightarrow B$ be a finite monomorphism of quasi-affinoid algebras and $A'$ be a quasi-affinoid $A$ algebra. Then $\phi$ induces a finite map $\psi:A'\rightarrow A'\otimes_A^s B$ and the nilradical of $A'$ contains Ker $(\psi)$.
\end{lem}
\begin{proof}

The first assertion is due to Proposition 5.4.8 of \cite{RoSPS}. For the second assertion, notice that because $\psi:A'\rightarrow A'\otimes_A^s B$ is a finite map extending $\phi$, by Theorem 9.3 of \cite{matsu} for each maximal ideal $\mfm$ of $A'$ there is a maximal ideal $\bar{\mfm}$ of $A'\otimes_A^s B$ such that $\bar{\mfm}\cap \psi(A')=\psi(\mfm)$. Let $a\in A'$ and assume there is a maximal ideal $\mfn$ of $A'$ such that $a\not\in \mfn$. Note that by the Nullstellensatz (Theorem 4.1.1 of \cite{RoSPS}), the nilradical of a quasi-affinoid algebra ring coincides with its Jacobson radical and therefore in this case  $a$ is not in the nilradical of $A'$. Then clearly $\psi(a)\not\in\bar{\mfn}$, and hence $a\not\in$ Ker($\psi$). 
\end{proof}

\begin{lem}\label{nil}
Let $A$ be a quasi-affinoid algebra and $c$ be a nilpotent element of $A\ls y_1,...,y_M\rs\lb\lambda_1,...,\lambda_N\rb$, then $c=\sum_{(\alpha,\beta)}a_{\alpha,\beta}y^{\alpha}\lambda^{\beta}$ where each $a_{\alpha,\beta}$ is a nilpotent element of $A$. Therefore if 
$$\phi:A\ls y_1,...,y_M\rs\lb\lambda_1,...,\lambda_N\rb\hookrightarrow B,$$
is a finite injection of quasi-affinoid algebras then $\phi$ induces a finite injection
$$\phi':A/\mathcal{N}(A)\ls y_1,...,y_M\rs\lb\lambda_1,...,\lambda_N\rb\hookrightarrow B/\mathcal{N}(B),$$
where $\mathcal{N}(A)$ and $\mathcal{N}(B)$ denote the nilradicals of $A$ and $B$ respectively.
\end{lem}
\begin{proof}
For the first assertion, first assume that $M=1$, $N=0$ or $M=0$, $N=1$. By assumption there is an $s$ such that $c^s=0$ and by Lemma \ref{idealmove} this implies that $a_0^s=0$ which in turn implies that $a_0\in\mathcal{N}(A)$. Therefore $c-a_0$ is also nilpotent and proceeding inductively we see that each $a_i$ is a nilpotent element of $A$. Now for the general case of arbitrary $M$ and $N$ we once again make use of Lemma \ref{idealmove} and argue inductively.

For the second assertion, finiteness of $\phi'$ is clear as $\mathcal{N}(A)\subset\phi^{-1}(\mathcal{N}(B))$. Let $c\in\phi^{-1}(\mathcal{N}(B))$ then $c$ is itself nilpotent and by the first assertion $c\in  \mathcal{N}(A)A\ls y_1,...,y_M\rs\lb\lambda_1,...,\lambda_N\rb$. Now the second assertion follows from Lemma \ref{idealmove}.
\end{proof}

\section{Restricted Dimension}

In this section we will discuss a restricted notion of Krull Dimension which is more closely related to the geometric properties of $D$-semianalytic sets than the Krull dimension (denoted {\em k-dim}).

\begin{defn}\label{Adim}
 Let $A$ be a quasi-affinoid algebra and $B$ a generalized ring of fractions over $A$ such that Dom$_AB\neq\emptyset$. For an ideal $I$ of $B$ we define the {\em Dom $A$-dimension of $I$ in $B$}, denoted $\text{dim}_A B/I$, as 
$$\text{dim}_A B/I:=\sup\{\text{k-dim }B_{\mathfrak{m}}/IB_{\mfm}: \mathfrak{m}\in\text{Dom}_A B\}$$
 If Dom$_AB=\emptyset$, then we define $\text{dim}_A B/I:=-1$. 
\end{defn}

\begin{rem} \label{Adimr}

i) Dom $A$-dimension of $I$ in $B$ behaves well under the operation of forming generalized rings of fractions. That is, if $B'$ is a generalized ring of fractions over $B$ (hence over $A$ also), then dim$_AB/I\geq$ dim$_AB'/IB'$. We prove this fact in Lemma \ref{dimeq}. 

ii) Notice that in these terms we can write 
$$\text{k-dim }B/I=\text{ dim}_B B/I.$$
so that the following results about Dom $A$-dimension have their analogues for Krull dimension of $B/I$.

iii) Notice also that $\text{dim}_{\smn} B/I=-1$ means that there is no maximal ideal $\mfm\in\text{Dom}_{\smn}B$ which contains $I$ and therefore in that case we have $\kdom B\cap V(I)_{K'}=\emptyset$ where $K'$ is a complete field extension of $K$.  On the other hand if  $\text{dim}_{\smn} B/I=0$, and $\mfp_1,..., \mfp_k$ are minimal primes of $I$ in $B$, then each $\mfp_i$ is either a maximal ideal of $B$ or there is no maximal ideal $\mathfrak{m}\supset\mfp_i$ such that $\mathfrak{m}\in\text{Dom}_{\smn}B$, hence $V(I)_{K'}\cap\kdom B$ is a finite or an empty set.
\end{rem}

Next we will establish basic facts about this restricted dimension to help us understand its connection with the geometry of the $D$-semianalytic sets.

\begin{lem}\label{minprim}
Let $A$ be a quasi-affinoid algebra and $B$ be a generalized ring of fractions over $A$, then
\begin{multline*}\text{dim}_A B/I=\sup \{\text{k-dim }B/\mfp: \text{ where } I\subset \mfp\in\text{Spec }B \text{ and there exists  }\\
\mathfrak{m}\in\text{Dom}_AB\text{ containing }\mfp\}.
\end{multline*}
\end{lem}
\begin{proof}
Write $A=\smn/J$ and $B=\Smn/J'$ and let $\bar{I}$ be the ideal corresponding to $I$ in $\Smn$. Assume $\mfp_0\supset \bar{I}\supset J'$ is a prime ideal and $\mathfrak{m}$ is a maximal ideal of $\Smn$ such that $\mathfrak{m}\in\text{Dom}_AB$ and such that $\text{dim}_A B/I=\text{k-dim }(\Smn)_{\mathfrak{m}}/\bar{I}_{\mfm}$. Now take a maximal prime ideal chain
$$\bar{I}_{\mfm}\subset\mfq_0\subset...\subset\mfq_d=\mfm_{\mfm}\subset(\Smn)_{\mfm}$$
such that $\text{dim}_A B/I=d$ and let $\mfp_0$ be the minimal prime divisor of $\bar{I}$ such that $(\mfp_0)_{\mfm}=\mfq_0$. Now assume that $\mathfrak{n}\subset\Smn$ is another maximal ideal lying over $\mfp_0$. By Lemma \ref{catenary} and Corollary 4.2.2 of \cite{RoSPS}, 
\begin{endproofeqnarray*}
\text{ht }\mfp_0&=&m+n+M+N-\text{k-dim }(\Smn)_{\mathfrak{m}}/(\mfp_0)_{\mfm}\\
&=&m+n+M+N-\text{k-dim }(\Smn)_{\mathfrak{n}}/(\mfp_0)_{\mfn}.
\end{endproofeqnarray*}
\end{proof}

Given an ideal $I$ of $\smn(E,K)$ it is natural to expect k-dim $\smn(E,K)/I$ to be the same as $\smn(E,K')/I\smn(E,K')$ where $K'$ is a complete field extension of $K$. In fact we can obtain a more general result as follows. Let $E'\subset {K'}^{\circ}$ be a complete, quasi-Noetherian ring and let $A$ be $\smn(E,K)/J$ for some ideal $J$ and assume 
$$B=A\ls y\rs\lb\lambda\rb/(\{g_iy_i-f_i\},\{G_j\lambda_j-F_j\})=\Smn(E,K)/J'$$ a generalized ring of fractions over $A$. Let $A'$ and $B'$ be results of ground field extension from $(E,K)$ to $(E',K')$ of $A$ and $B$ respectively. That is
$$A':=\smn(E',K')/J\cdot\smn(E',K'),$$ 
and 
$$B'=\Smn(E',K')/J'\cdot\Smn(E',K').$$ 
By Lemma 4.2.8 of \cite{RoSPS}, $\smn(E',K')$ is faithfully flat over $\smn(E,K)$, therefore if $\mfm\in$ Dom$_AB$, then there is an $\mfm'\in$ Dom$_{A'}B'$ such that $\mfm'\cap B=\mfm$ and we can state the following corollary. 

\begin{cor}\label{extension}
$\text{dim}_AB/I=\text{dim}_{A'}B'/IB'$.
\end{cor}
\begin{proof}
Observe that the following is a consequence of \cite{matsu}, Theorem 15.1.

{\em Claim.} Suppose $\mfm\in\text{Dom}_{A'}B'$ and $\mfm\supset IB'$. Put $\mfn:=B\cap\mfm$ and suppose $\mfn$ is a maximal ideal. Then $\mfn\in\text{Dom}_AB$, $\mfn\supset I$ and 
$$\text{k-dim }(B/I)_{\mfn}=\text{k-dim }(B'/IB')_{\mfm}.$$

In particular, since the map $B\rightarrow B'$ is faithfully flat, from the claim and the discussion above we have
$$\text{dim}_AB/I\leq\text{dim}_{A'}B'/IB'.$$
On the other hand, by Lemma \ref{minprim} we may assume that 
$$\bar{I}\cdot\Smn(E',K')\subset\mfq_0\subset...\subset\mfq_d\subset\Smn(E',K')$$
is a maximal chain of prime ideals with $d=\text{dim}_{A'}B'/IB'$ and $\mfq_d\in\text{Dom}_{A'}B'$. Thus by the Going Down Theorem and faithful flatness, $\mfp_0:=\mfq_0\cap\smn$ is a minimal prime divisor of $\bar{I}$.

Let $\mfq$ be a minimal prime divisor of $\bar{I}\cdot\Smn(E',K')$. By \cite{matsu}, Theorem 15.1
$$\begin{array}{rcl}
\text{ht }\mfq&=&\text{ht }\mfp_0+\text{dim }(\Smn(E',K'))_{\mfq}/\mfp_0(\Smn(E',K'))_{\mfq}\\
&=&\text{ht }\mfp_0
\end{array}$$
since $\mfq$ is a minimal prime divisor of $\bar{I}\cdot\Smn(E',K')$, hence also of $\mfp_0\cdot\Smn(E',K')$. So $d=\text{dim }\Smn/\mfq$ for any such $\mfq$. Thus by the {\em Claim} above and Lemma \ref{minprim}, it suffices to show that there is a minimal prime divisor $\mfq$ of $\mfp_0\Smn(E',K')$ such that $\mfq\subset\mfm$ for some $\mfm\in\text{Dom}_{A'}B'$ such that $\mfn:=\mfm\cap B$ is a maximal ideal. Since $\text{Dom}_AB$ is a Zariski-open set containing a point on $\mfp_0$, this follows from the faithful flatness.
\end{proof}

In the rest of this section we will investigate the behavior of the restricted dimension under the operation of forming generalized rings of fractions. By Lemma \ref{idealmove}, the following is easily proved.

\begin{lem}\label{+1}
Let $A=\smn/J$ be a quasi affinoid algebra, $I\subset A$ a proper ideal, $z$ a variable not appearing in $\smn$. Then 
$$\text{k-dim }A\ls z\rs/IA\ls z\rs=\text{ k-dim }A\lb z\rb/IA\lb z\rb=\text{ k-dim }A/I+1.$$
\end{lem}

\begin{lem}\label{dimeq}
Let $I$ be an ideal of $\smn$, $f,g$ two elements of $\smn$ and $\mathfrak{m}$ be a maximal ideal of $\smn\ls z\rs$ containing $I$ and $gz-f$ but not $g$. Let $\mathfrak{n}=\smn\cap\mfm$, then k-dim $(\smn)_{\mfn}/I(\smn)_{\mfn}=$ k-dim $(\smn\ls z\rs)_{\mfm}/(I\cup\{gz-f\})(\smn\ls z\rs)_{\mfm}.$ The same statement is true if we replace $\smn\ls z\rs$ with $\smn\lb z\rb$.
\end{lem}
\begin{proof}
Once again we will only give the proof for  $\smn\ls z\rs$. It is easy to see that 
$$\text{k-dim }(\smn\ls z\rs)_{\mfm}/I(\smn\ls z\rs)_{\mfm}=\text{k-dim }(\smn)_{\mfn}/I(\smn)_{\mfn}+1.$$ 
To show that 
$$\text{k-dim }(\smn)_{\mfn}/I(\smn)_{\mfn}\geq\text{k-dim }(\smn\ls z\rs)_{\mfm}/(I\cup\{gz-f\})(\smn\ls z\rs)_{\mfm},$$ 
it is enough to show that given a maximal chain
$$(I\cup\{gz-f\})\cdot\smn\ls z\rs\subset\mfp_0\subset...\subset\mfp_d=\mfm,$$
we can construct a strict chain 
$$I\smn\ls z\rs\subset\mfq_0\subset...\subset\mfq_{d+1}=\mfm.$$

Notice that no minimal prime of $I \smn\ls z\rs$ which is contained in $\mfm$ contains $gz-f$. This is because by Lemma \ref{idealmove} such minimal primes are of the form $\mfq \smn\ls z \rs$ for some minimal prime $\mfq$ of $I$ in $\smn$. Indeed if it is the case that $b_1\sum a_{1,i}z^i+...+b_k\sum a_{k,i}z^i=gz-f$ with $b_j\in\mfq$, $a_{j,i}\in \smn$, then we have $g=b_1a_{1,i}+...+b_ka_{k,i}\in\mfq$ contradicting $g\not\in\mfm$. Now let $\mfq_0$ be a minimal prime divisor of $I\smn\ls z\rs$ which is contained in $\mfp_0$. By the explanation above paragraph, $gz-f\not\in\mfq_0$, and therefore $\mfq_0\neq\mfp_0$. Hence setting $\mfq_i:=\mfp_{i-1}$ for $1\leq i\leq d+1$ gives us the desired prime ideal chain. The argument for $\smn\lb z\rb$ is similar.
\end{proof}

 In other words, the above lemma states that if dim$_A B/I=$ k-dim $(B/I)_{\mfm}=d$ for some maximal ideal $\mfm\supset I$ in Dom$_A B$, and $g\not\in\mfm$, then $$\text{dim}_A B/I=\text{ dim}_A B\ls f/g\rs/IB\ls f/g\rs,$$ 
and of course similarly for $B\lb f/g\rb$. As an immediate consequence we also see that dim$_AB/I\leq$ k-dim $A$.

\section{Manifolds and Geometric Dimension}

In this section we will investigate the geometric properties of $D$-semianalytic sets more closely. Although we aim to get smooth pieces of $D$-semianalytic sets which are manifolds, we will not try to make the concept of manifolds precise to save time as all the $K$-$n$-manifolds we consider will be $D$-semianalytic sets which are locally graphs of $n$-dimensional open balls $U$ in $(K^{\circ})^n$ under power series converging on $U$ with coefficients from $K$. The reader can find a more general treatment of $K$-analytic manifolds in Subsection 1.4 of \cite{strat}. Working with such manifolds, our main tool will be the next well known theorem. For a proof we refer the reader to Proposition 10.8 of \cite{abh}. 

\begin{thm}[Implicit Function Theorem]\label{implicit}
Let $f_1,...,f_n$ be elements of the Tate Algebra $T_m=K\ls x_1,...,x_m\rs$, $n\leq m$ and $\Delta$ be the determinant of the $n$-dimensional minor $\partial (f_1,...,f_n)/\partial (x_1,...,x_n)$ of the Jacobian Matrix $\partial f/\partial x $. Let $\pi$ be the projection map onto the space corresponding to variables $x_{n+1},...,x_m$. Then for all $\bar{p}\in V(f_1,...,f_n)_K$ satisfying $ \Delta(\bar{p})\neq 0$ there is a neighborhood $U$ of $\bar{p}$ such that $U\cap V(f_1,...,f_n)_K$ is the graph of a function on $\pi(U)$ given by power series convergent on $\pi(U)$ with coefficients from $K$.
\end{thm}

This theorem is especially useful when it is used in conjunction with the following slightly modified results (Theorem 3.1 and 3.2) from \cite{strat}.

\begin{lem}\label{nicegen0}
Suppose {\em Char }$K=0$ and $B$ is a generalized ring of fractions over $\smn$, Let $f_1,...,f_k\in B$, $I$ be the ideal generated by $f_1,...,f_k$, and let $\mfp$ be a minimal prime divisor of $I$ such that there is some $\mfn\in\text{Dom}_{\smn} B$, $\mfn\supset\mfp$, {\em ht} $\mfp=r$, then there are differential polynomials $P_1,...,P_r,Q$ with integer coefficients  and a maximal ideal $\mfm\in\text{Dom}_{\smn} B$ such that

i) $Q(f)\not\in\mfm$ and $P_1(f)/Q(f),...,P_r(f)/Q(f)$ generate $\mfp B_{\mfm}$.

ii) some $r\times r$ minor of the Jacobian matrix $\frac{\partial}{\partial x,\rho}(P_1(f)/Q(f),...,P_r(f)/Q(f))$ does not belong to $\mfp B_{\mfm}$. 
\end{lem}

\begin{lem}\label{nicegenp}
Suppose {\em Char }$K=p>0$ and $B$ is a generalized ring of fractions over $\smn$, Let $f_1,...,f_k\in B$, $I$ be the ideal generated by $f_1,...,f_k$, and let $\mfp$ be a minimal prime divisor of $I$ of height $r$ such that there is some $\pb\in\dom B$, let $\mfm$ be the corresponding maximal ideal of $B$. Then there is an $l\in\mathbb{N}$, an $R$-domain $U$ such that $\pb$ is an element of the set of $K$-rational points in $U$ and there exist Hasse differential polynomials $P_1,...P_r,Q$, with coefficients in $\mathbb{F}_p$, such that
$$Q(f),P_i(f)\in\mathcal{O}_{K^{1/p^l}}(U)^{p^l},$$
where $\mathcal{O}_{K^{1/p^l}}(U):=(E,K^{1/p^l})\otimes^s_{(E,K)}\mathcal{O}(U)_K$, and

i) $Q(f)\not\in\mfm$ and the nilradical of the ideal of $B_{\mfm}$ generated by $P_1(f)/Q(f),...,$ $P_r(f)/Q(f)$ is $\mfp B_{\mfm}$, and

ii) Some $r\times r$ minor of the Jacobian $\frac{\partial}{\partial x,\rho}(P_1(f)/Q(f),...,P_r(f)/Q(f))^{1/p^l}$ does not belong to $\mfq$ for any prime ideal $\mfq$ of $\mathcal{O}_{K^{1/p^l}}(U)$ lying above $\mfp$.
\end{lem}

Next we will define a geometric notion of dimension to continue our investigation.
\begin{defn}\label{gdim}
 Define the {\em geometric dimension}, g-dim $X$,  of a nonempty set $X\subset K^m$ to be the greatest integer $d$ such that the image of $X$ under coordinate projection onto a $d$ dimensional coordinate hyper-plane has an interior point. For $X=\emptyset$ define g-dim $X=-1$. 

For $\bar{p}\in K^m$ define g-dim $X_{\pb}$ to be the minimum of g-dim $X\cap U$ where $U$ is a neighborhood of $\pb$.
\end{defn}

Let $B$ be a generalized ring of fractions over $\smn$, and let $I$ be an ideal of $B$. We use the customary notation $\mathcal{I}(V(I)_{K'}\cap\kdom B)$ to denote the ideal $\{f\in B: f(\pb)=0 \text{ for all }\pb\in V(I)_{K'}\cap\kdom B\}$ of $B$. In the rest of this paper we establish the connections between the geometric and the weak dimensions (Definition \ref{wdim}) of $\dom B\cap V(I)_K$ and the restricted dimension dim$_{\smn} B/I$. As the restricted dimension is defined in terms of the Krull dimension, it is more closely related to the geometric and weak dimensions of $\bar{K}\text{-Dom}_{m,n}B\cap V(I)_{\bar{K}}$ where $\bar{K}$ is an algebraically closed complete extension of $K$. In order for the comparison between these different types of dimension to be meaningful, we need to assume $I=\mathcal{I}(V(I)_K\cap\dom B)$ so that $I$ is the largest among the ideals of $B$ which define the same $D$-semianalytic set $\dom B\cap V(I)_K$. If that is the case for an ideal $I$, then it is also the case for the minimal prime divisors of $I$ as the next lemma shows.

\begin{lem}\label{primdivmax}
Suppose $B$ is a generalized ring of fractions over $\smn$, and $I$ is an ideal of $B$ such that $I=\mathcal{I}(V(I)_{K'}\cap\kdom B)$. Then $I$ is a radical ideal and if $I=\cap_{i=1}^s\mfp_i$ is an irredundant minimal prime decomposition of $I$, then $\mfp_i=\mathcal{I}(V(\mfp_i)_{K'}\cap\kdom B)$.
\end{lem}
\begin{proof}
That $\mathcal{I}(V(I)_{K'}\cap\kdom B)$ is a radical ideal is clear. The rest of the assertion is also clear if $I$ is prime, hence assume $s>1$. Write $J_i:=\mathcal{I}(V(\mfp_i)_{K'}\cap\kdom(B))$, and assume $f_i\in J_i\backslash\mfp_i$. As we assumed $\cap_{i=1}^s\mfp_i$ is an irredundant minimal prime decomposition of $I$, we can find $g_i$ which is an element of  $(\bigcap_{j\neq i}\mfp_j)\backslash\mfp_i$. Observe that for all $\pb\in V(I)_{K'}\cap\kdom B$, $g_i(\pb)\neq 0$ implies $f_i(\pb)=0$, hence $f_ig_i$ is an element of  $\mathcal{I}(V(I)_{K'}\cap\kdom B)$ which is equal to $I$ by assumption. Thus $f_ig_i\in\mfp_i$, a contradiction.
\end{proof}

As a corollary to the above lemma we see that if $I=\mathcal{I}(V(I)_{K'}\cap\kdom B)$ then for each minimal prime ideal $\mfp$ of $I$ there is a maximal ideal $\mfm\in\text{ Dom}_{\smn}B$ containing $\mfp$. Therefore k-dim $B/I=$ dim$_{\smn} B/I$.

The Smooth Stratification Theorem for subanalytic subsets of $(\bar{K}^{\circ})^m\times (\bar{K}^{\circ\circ})^n$ is proved (Theorem 4.4) in \cite{strat} where $\bar{K}$ is an algebraically closed complete extension of $K$. We would like to remind the reader that those sets are also $D$-semianalytic by the Quantifier Elimination Theorem (Theorem 3.8.1) of \cite{rigid}. The key step in proving the Smooth Stratification Theorem is Corollary 3.3 of \cite{strat} and we are going to generalize this corollary to $D$-semianalytic subsets of $(K^{\circ})^m\times (K^{\circ\circ})^n$.

\begin{thm}[Smooth Stratification of $D$-semianalytic Sets]\label{stratify}
Assume {\em Char }$K=0$. Let $B$ be a generalized ring of fractions over $\smn$, $I\subset B$ an ideal which satisfies $I=\mathcal{I}(V(I)_K\cap\dom  B)$, then there are $D$-semianalytic $K$-manifolds $X_1,...,X_r$ (not necessarily disjoint) such that 
 $$V(I)_K\cap\dom B=X_1\cup...\cup X_r$$
and
 $\max_i$ g-dim $X_i=$ dim$_{\smn}B/I$= k-dim $B/I$.
\end{thm}
\begin{proof}
We will actually only prove that $\max_i$ g-dim $X_i\geq$ dim$_{\smn}B/I$ and the other inequality will follow from Lemma \ref{rdimdomwdim}.

The assertion is clear if there are no points in $\dom B\cap V(I)_K$ in which case $I=(1)$. So assume there is a $\pb\in\dom B\cap V(I)_K$. By Lemma \ref{primdivmax} each minimal prime divisor $\mfp_i$ of $I$ satisfies $\mfp_i= \mathcal{I}(V(\mfp_i)_K\cap \dom B)$, hence we may assume that $I$ is prime. Let $r$ be {\em ht }$I$. We will proceed by induction on dim$_{\smn}B/I$.

Write $f=(f_1,...,f_k)$ for generators of $I$ and apply Lemma \ref{nicegen0} to $I$ to get a maximal ideal $\mfm\in\text{Dom}_{\smn} B$, $I\subset \mfm$, together with differential polynomials $P_1,...,P_r,Q$ with integer coefficients and the determinant of an $r\times r$ minor $\Delta$ of the Jacobian $\frac{\partial}{\partial (x,\rho)}(P_1(f)/Q(f),...,P_r(f)/Q(f))$, such that $\Delta\not\in IB_{\mfm}$. Notice that because $I$ is prime $P_i(f)\in I$ for all $i$. In order to simplify the notation we may assume that 
\begin{displaymath}\Delta=\left|\begin{array}{ccc}
\frac{\displaystyle \frac{\partial P_1(f)}{\partial x_1}Q(f)-\frac{\partial Q(f)}{\partial x_1}P_1(f)}{\displaystyle Q(f)^2}&\cdots&\frac{\displaystyle \frac{\partial P_{r}(f)}{\partial x_1}Q(f)-\frac{\partial Q(f)}{\partial x_1}P_r(f)}{\displaystyle Q(f)^2}\\
\\
\vdots&\ddots&\vdots\\
\\
\frac{\displaystyle \frac{\partial P_1(f)}{\partial x_r}Q(f)-\frac{\partial Q(f)}{\partial x_r}P_1(f)}{\displaystyle Q(f)^2}&\cdots&\frac{\displaystyle \frac{\partial P_{r}(f)}{\partial x_r}Q(f)-\frac{\partial Q(f)}{\partial x_r}P_r(f)}{\displaystyle Q(f)^2}
\end{array}\right|.\end{displaymath}

Observe that each $\frac{\partial Q(f)}{\partial x_i}P_j(f)/Q(f)^2$ is in $IB_{\mfm}$. Therefore by elementary properties of determinant, $\Delta= \Delta'/Q(f)^r+g$ where $g\in IB_{\mfm}$ and
\begin{displaymath} \Delta'=\left|\begin{array}{ccc}
\displaystyle\frac{\partial P_1(f)}{\partial x_1}&\cdots&\displaystyle\frac{\partial P_{r}(f)}{\partial x_1}\\
\\
\vdots&\ddots&\vdots\\
\\
\displaystyle\frac{\partial P_1(f)}{\partial x_r}&\cdots&\displaystyle\frac{\partial P_{r}(f)}{\partial x_r}
\end{array}\right|.\end{displaymath}
Therefore $\Delta\not\in IB_{\mfm}$ implies $\Delta'\not\in I$.  
 
Now define 
$$X:=\{\bar{q}\in \dom B\cap V(I)_K: \Delta'(\bar{q})\neq 0\}.$$ 
For all $\qb\in X$ and maximal ideal $\mfn$ corresponding to $\qb$ we have, by Theorem 30.4 of \cite{matsu}, $P_1(f),...,P_r(f)$ generate $I B_{\mfn}$ (it is easy to show that for an $\mfn\in\text{Dom}_{\smn} B$, $B_{\mfn}$ is regular). Notice that, in this case there is a rational neighborhood $W$ of $\qb$ such that  $P_1(f),...,P_r(f)$ generate $I\mathcal{O}(W)$ and by Theorem \ref{implicit}, there is a rational neighborhood $U\subset W$ of $\qb$ such that $V(I)_K\cap U$ is the graph of analytic functions over some $m+n-r$ dimensional open poly-disc. Notice also that $X\neq\emptyset$ by the assumption that $I=\mathcal{I}(V(I)_K\cap\dom B)$ and  because $\Delta'\not\in I$.

On the other hand $(V(I)_K\cap\dom B)\backslash X=V(I\cup\{\Delta'\})_K\cap\dom B$, and because $I$ is prime, we see that dim$_{\smn}B/(I\cup\{\Delta'\})<$ dim$_{\smn}B/I$, and the result follows by induction.
\end{proof}

Note that our main tool in the above proof, Lemma \ref{nicegen0} has an analogue, Lemma \ref{nicegenp}  for the case {\em Char $K=p>0$}. However Lemma \ref{nicegenp} is much weaker than Lemma \ref{nicegen0} in the sense that the generators $P_i(f)/Q(f)$ we get have coefficients in some $K^{1/p^l}$ and they work only locally. Therefore an analogue of Theorem \ref{stratify} can be proved similarly but with the additional condition $[K:K^p]<\infty$ and with a local, rather than a global conclusion.

\begin{thm}\label{stratifyp}
Assume {\em Char }$K=p>0$ and $[K:K^p]<\infty$. Let $B$ be a generalized ring of fractions over $\smn$, $I\subset B$ an ideal of $B$, then for all $\pb\in \dom B\cap V(I)_K$, there is an $R$-domain $U$ such that $\pb$ is an element of the set $W$ of $K$-rational points in $U$ and such that there exist $D$-semianalytic $K$-manifolds $X_1,...,X_r$ (not necessarily disjoint) such that 
 $$V(I)_K\cap W=X_1\cup...\cup X_r.$$
Furthermore if $J=\mathcal{I}(W\cap V(I)_K)$ is an ideal of $\mathcal{O}(U)_K$, then
 $\sup_i$ g-dim $X_i\geq$ k-dim $\mathcal{O}(U)_K/J$.
\end{thm}

The geometric dimension we have been using so far has the advantage that it is based on open sets and the intersection of open sets with hyperplanes are open subsets of those hyperplanes. This observation often helps us in determining the dimension of the intersection of a $D$-semianalytic set with a hyperplane. On the other hand, one naturally expects that the dimension of a union of a finite number of sets is the maximum of the dimensions of the individual sets. For this purpose it is convenient to use the next definition of dimension.

\begin{defn}\label{wdim}
For a nonempty subset $X$ of $K^m$ define the {\em weak dimension}, w-dim $X$, to be the greatest integer $d$ such that the image of $X$ under coordinate projection onto a $d$ dimensional coordinate hyper-plane is somewhere dense. If $X$ is empty, define w-dim $X=-1$. For $\pb\in K^m$, w-dim $X_{\pb}$ is defined the same way as in geometric dimension.
\end{defn}

It is clear that the weak dimension of a set dominates the geometric dimension. In Section 6 we explore further relationships between the geometric and the weak dimension and prove that in characteristic $0$, weak dimension of a $D$-semianalytic set is equal to the geometric dimension of the set (Theorem \ref{alleq0}). For the moment there is an easy yet useful lemma that we can obtain right away.

\begin{lem}\label{wdimgrdimeq1}
Let $X=\dom B\cap V(I)_K$, for a generalized ring of fractions $B$ over $\smn$ and $I$ an ideal of $B$ satisfying $I=\mathcal{I}(\dom B\cap V(I)_K)$. If w-dim $X=m+n$, then g-dim $X=$ dim$_{\smn}B/I=m+n$.
\end{lem}
\begin{proof}
If w-dim $X=m+n$, then since $X$ is a locally closed somewhere dense subset of $(K^\circ)^m\times(K^{\circ\circ})^n$, $X$ has to contain an interior point $\pb$. We may assume that $\pb$ is the origin and the open ball with radius $|\varepsilon|$ around the origin is contained in $X$ for some $\varepsilon\in K^\circ$. Notice that 
$$T_{m+n}\simeq B':=B\ls x_1/\varepsilon,...,x_m/\varepsilon,\rho_1/\varepsilon,...,\rho_n/\varepsilon\rs$$ 
and by Lemma \ref{dimeq}, dim$_{\smn}B/I\geq$ k-dim $B'/IB'$. Now the result follows from the fact that if $J\subset T_{m+n}$ and $V(J)_K=(K^\circ)^{m+n}$ then $J=(0)$.
\end{proof}

\section{Normalization}

In this section we obtain a normalization lemma that enables us to compare the geometric dimension of a $D$-semianalytic set with the restricted dimension of the associated quasi-affinoid algebras, which we do in Section 6. Note that there is no ``global'' normalization for quasi-affinoid algebras in general. That is, given a quasi-affinoid algebra $A=\smn/J$, it is not always possible to find integers $m'$, $n'$ and a finite monomorphism $\phi:S_{m',n'}\rightarrow A$. An example of failure of normalization can be found in Example 2.3.5 of \cite{RoSPS} where it was shown that it is not possible to normalize the quasi-affinoid algebra $A=S_{1,1}/(x\rho)$ in the sense above. Nevertheless in Lemma \ref{norm1} we show that given a quasi-affinoid algebra $B$ it is possible to obtain finitely many normalized (in a slightly different sense than the above one) quasi-affinoid algebras $B_j$ with the union of the corresponding $D$-semianalytic sets (see Definition \ref{Dsemi}) equal to the $D$-semianalytic set corresponding to $B$, and the finite monomorphisms $\phi_j$ act as Weierstrass changes of variables. On the other hand one may have a quasi-affinoid algebra $A$, in whose presentation some of the variables are designated as parameter variables. For example this will be the case for us when we prove Theorem \ref{alleq0} and \ref{dimadd} on dimensions of $D$-semianalytic sets as we will make such a distinction between the variables corresponding the target space of a projection map $\pi$ and the variables corresponding the fiber space. The maps $\phi_j$ we find in Lemma \ref{norm1} will also fix all the parameter variables and hence we call Lemma \ref{norm1} the Parameterized Normalization Lemma. This lemma  will also play an essential part in the proof of a parameterized version of Theorem \ref{stratify} in a subsequent paper \cite{parastrat}, in which we will generalize the main theorem of \cite{bart} by Bartenwerfer. The proof of Lemma \ref{norm1} closely follows the ideas in the proof of Theorem 6.1.2 (Finiteness Theorem) of \cite{RoSPS}.

 Although we name Lemma \ref{norm1} as the Parameterized Normalization Lemma,  all of the statements and observations in this section are about refinements of the process of normalization. Our results are distributed in several separate statements for ease in reading. Hence, when we refer to the Parameterized Normalization Lemma as opposed to referring to Lemma \ref{norm1} elsewhere in this paper, we trust that the reader will take this to be a referral to all of the statements and observations in this section. 

 Before we state the main result of this section, we will generalize $R$-domains (see Definition \ref{rdom}) to the parameterized case. Although these domains are not essential for our discussion about dimension in general, because $R$-domains play an important role in non-Archimedean geometry and because our results are expressed in full generality in terms of these domains, we will take time to define them. 

First let us observe a nice property about the $K$-rational points in a domain of generalized ring of fractions which helps justify the terminology {\em parameterized families of $R$-domains} below. Let $A$ be a generalized ring of fractions over $\smn$ and $f\in A\ls y_1,...,y_M\rs\lb\lambda_1,...,\lambda_N\rb$, then we can write $f=\sum_{\alpha,\beta} a_{\alpha\beta}y^{\alpha}\lambda^{\beta}$ where $a_{\alpha,\beta}\in A$. Let $\pb\in\text{ Dom}_{\smn}A$, then 
$$f(\pb,y,\lambda):=\sum_{\alpha,\beta}a_{\alpha,\beta}(\pb)y^{\alpha}\lambda^{\beta}$$ 
is an element of $S_{M,N}$. This fact easily follows from Definition \ref{separ}.

\begin{defn}\label{PFOR}

Let $S_{m+M,n+N}$ denote the ring of separated power series over the variables $x_1,...,x_m;$ $y_1,...,y_M;$ $\rho_1,...,\rho_n$ and $\lambda_1,...,\lambda_N$. Define a {\em parameterized family of R-domains} ($PRD$) {\em over $\Smn$ with parameters $x_1,...,x_m$ and $\rho_1,...,\rho_n$}  defined by {\em the parameter ring $A$ and the function ring $B$} (over $(E,K)$, see Definition \ref{separ}) inductively as follows:

 i) Max $\Smn$ is a $PRD$ over $\Smn$ defined by function ring $B=\Smn$ and parameter ring $A=\smn$.

 ii) Suppose $X$ is a $PRD$ over $\Smn$ defined by parameter ring $A$, and function ring $B$, and let $A'$ be a generalized ring of fractions over $A$, suppose that $(f_1,...,f_k,h_1,...,h_l,g)$ generate the  unit ideal of $B$, put 
$$B':=A' \otimes_A^s B\ls f_1/g,...,f_k/g \rs\lb h_1/g,...,h_l/g \rb.$$
Then $X'=\text{Dom}_{\Smn}B'$ is a $PRD$ over $\Smn$ defined by the parameter ring $A'$, and the function ring $B'$.
\end{defn}

In other words, a $PRD$ is obtained by relaxing the conditions in the inductive construction of generalized rings of fractions whose domains are $R$-domains. This is in the sense that we allow adjoining arbitrary fractions from the parameter ring $A$. The resulting domains are not necessarily $R$-domains, but their specializations at suitable points from the parameter space are, as the following remark indicates.

\begin{rem}\label{PFORr} 
If $X$ is a $PRD$ over $\Smn$ defined by the parameter ring $A$, and the function ring $B$ then by Lemma \ref{idealmove} we can write
\begin{multline}
\label{pfr}\tag{$\dag$}
B=A\ls y_1,...,y_{M+S}\rs\lb\lambda_1,...,\lambda_{N+T}\rb/(\{g_{M+i}y_{M+i}-h_{M+i}\}_{i=1}^S\cup\\ 
\{G_{N+j}\lambda_{N+j}-H_{N+j}\}_{j=1}^T).
\end{multline}

Let $\bar{p}\in\dom A$ and define 
\begin{multline*}
B_{\pb}=S_{M+S,N+T}/(\{g_{M+i}(\pb,y,\lambda)y_{M+i}-h_{M+i}(\pb,y,\lambda)\}_{i=1}^S\cup \\
\{G_{N+j}(\pb,y,\lambda)\lambda_{N+j}-H_{N+j}(\pb,y,\lambda)\}_{j=1}^T).
\end{multline*}
Then Dom$_{S_{M,N}}B_{\pb}$ is an $R$-domain over $S_{M,N}$. 
\end{rem}

The next lemma and its proof describes the main body of the normalization process.

\begin{lem}[Parameterized Normalization Lemma]\label{norm1}
Adopting the notation of Definition \ref{PFOR} let $X$ be a $PRD$ over $\Smn(E,K)$ defined by the function ring  $B$ and the parameter ring $A$. Write $B$ as in Equation \eqref{pfr} and let  $I$ be  an ideal of $B$, then there exist finitely many parameterized families of $R$-domains $X_j$ defined by the function rings $B_j$, parameter rings $A_j$ (which are generalized rings of fractions over $\Smn(E,K)$  and $\smn(E,K)$ respectively), together with ideals $I_j\subset B_j$ satisfying $IB_j\subset I_j$ and integers $M_j$, $N_j$ such that for any complete extension $K'$ of $K$

i) $\kDom B=\bigcup_{j}\kDom B_j$,

ii) $\kDom B\cap V(I)_{K'}=\bigcup_j (\kDom B_j\cap V(I_j)_{K'})$, 

iii) $\kDom B_j\cap V(I_j)_{K'}\cap \kDom B_i=\emptyset\text{ for }i\neq j,$ 

iv) writing each $B_j$ in the form 
\begin{multline*} B_j=A_j\ls y_1,....,y_{M+S_j}\rs\lb\lambda_1,...,\lambda_{N+T_j}\rb/(\{y_{M+k}g_{j,M+k}-h_{j,M+k}\}_{k=1}^{S_j}\cup\\
\{\lambda_{N+l} G_{j,M+l}-H_{j,M+l}\}_{l=1}^{T_j}),
\end{multline*}
after Weierstrass changes of variables among $y$'s and among $\lambda$'s separately, there are finite monomorphisms 
$$\phi_j:A_j/(I_j\cap A_j)\ls y_1,...,y_{M_j}\rs\lb\lambda_1,...,\lambda_{N_j}\rb\hookrightarrow B_j/I_j.$$
\end{lem}

\begin{proof} The proof will be in two steps and first we will assume $X$ is relatively simple.

{\em Case 1.} 
$$B=A\ls y_1,...,y_M\rs\lb \lambda_1,...,\lambda_N\rb.$$

We will prove this case by induction on the pair $(M,N)$. Assume the statement holds for all $B'=A'\ls y_1,...,y_{M'}\rs\lb \lambda_1,...,\lambda_{N'}\rb$ for some generalized ring of fractions $A'$ over $\smn$ with $(M',N')<(M,N)$ in the lexicographic order. Let $f_1,...f_k$ generate $I$ and write $f_i=\Sigma a_{i\mu\nu}(x,\rho)y^\mu \lambda^\nu$, where $a_{i\mu\nu}(x,\rho)\in A$. Let $\bar{I}$ be the ideal of $B$ generated by $\{a_{i\mu\nu}\}$. Then by Lemma \ref{flat} $\bar{I}\cap A$ is the ideal generated by the elements $\{a_{i\mu\nu}\}$ and by Lemma \ref{idealmove} $I\subset\bar{I}$ hence the assertion (iv) of the lemma is satisfied for $A_0=A$, $B_0=B$, $X_0=X$, $I_0=\bar{I}$, $M_0=M$, $N_0=N$ and the identity map $\phi_0$.
 
Next we will form finitely many generalized rings of fractions such that on each of them one of the $f_1,...,f_k$ will be preregular in the sense of Definition 2.3.7 of \cite{RoSPS} and domain of each of them will be a $PRD$. Notice that by Lemma 3.1.6 of \cite{RoSPS}, there is a finite index set $Z\subseteq \mathbb{N}^M\times \mathbb{N}^N$, such that for each $\bar{p}\in \kdom(A)\setminus V(\bar{I}\cap A)_{K'}$ there is an $i_0$, $1\leq i_0\leq k$ and an index $(\mu_0,\nu_0)\in Z$ such that
$$\begin{array}{rcl}
|a_{i_0\mu_0\nu_0}(\bar{p})|&\geq&|a_{i\mu\nu}(\bar{p})|\text{ for all }i,\mu,\nu\\
|a_{i_0\mu_0\nu_0}(\bar{p})|&>&|a_{i_0\mu\nu}(\bar{p})|\text{ for all }\nu<\nu_0\text{ and all }\mu\\
|a_{i_0\mu_0\nu_0}(\bar{p})|&>&|a_{i_0\mu\nu_0}(\bar{p})|\text{ for all }\mu>\mu_0.\\
\end{array}$$ 

For every $i$, $1\leq i \leq k$ and $(\mu,\nu)\in Z$, define
$$A_{i\mu\nu}=A\ls\left\{ \frac{a_{l\alpha\beta}}{a_{i\mu\nu}}\right\}_{l\alpha\beta} \rs\lb\left\{ \frac{a_{i\gamma\delta}}{a_{i\mu\nu}}\right\}_{\gamma\delta}\rb$$
where $(\gamma,\delta)$ runs through all pairs in $Z$ with $\delta<\nu$ or with $\alpha>\mu$ and where $(l,\alpha,\beta)$ runs through all remaining tuples in $\{1,...,k\}\times Z$. 

For all $(i,\mu,\nu)\in \{1,...,k\}\times Z$, also define 
$$B_{i\mu\nu}=A_{i\mu\nu}\ls y_1,...,y_M \rs\lb\lambda_1,...,\lambda_N\rb,$$
and put $X_{i\mu\nu}=\text{Dom}_{\Smn}B_{i\mu\nu}$. Notice that  each $X_{i\mu\nu}$ is a $PRD$ and by the discussion above we have 
\begin{multline*}\kDom(B)\cap V(I)_{K'}=(\kDom B\cap V(\bar{I})_{K'})\cup\\
\bigcup_{i\mu\nu}(\kDom B_{i\mu\nu}\cap V(IB_{i\mu\nu})_{K'}),
\end{multline*}
and  $\kDom B\cap V(\bar{I})_{K'}\cap \kDom B_{i\mu\nu}=\emptyset$.

Notice that although $a_{i\mu\nu}^{-1}$ may not be a member of $B_{i\mu\nu}$, again by Lemma 3.1.6 of \cite{RoSPS} (strong Noetherianness) we have $a_{i\mu\nu}^{-1}f_{i}\in B_{i\mu\nu}$ and it vanishes on every point in $\kdom B_{i\mu\nu}\cap V(I)_{K'}$ so that by replacing $I$ with $I\cup \{a_{i\mu\nu}^{-1}f_{i}\}$, we may assume that $a_{i\mu\nu}=1$. Notice also that  
$$\kDom B_{i\mu\nu}\cap \kDom B_{i\alpha\beta}=\emptyset,$$
for all $(\mu,\nu)$, $(\alpha,\beta)\in Z$, with $(\mu,\nu)\neq (\alpha,\beta)$ and in fact by further subdividing into $PRD$s we may assume that the $\kDom B_{i\mu\nu}$ are pairwise disjoint.

Next we will concentrate on the individual $PRD$s defined above. Fix an $i_0\mu_0\nu_0\in\{1,...,k\}\times Z$ and define
$$g_{i_0\mu_0\nu_0}=\underset{\mu}{\sum} a_{i_0\mu\nu_0}y^{\mu}.$$
Using $g_{i_0\mu_0\nu_0}$, we define two $PRD$s 
$$\begin{array}{rcl}
V_{i_0\mu_0\nu_0}&=&\text{Dom}_{\Smn}B_{i_0\mu_0\nu_0}\ls y_{M+1}\rs/(y_{M+1}g_{i_0\mu_0\nu_0}-1)\\
W_{i_0\mu_0\nu_0}&=&\text{Dom}_{\Smn}B_{i_0\mu_0\nu_0}\lb \lambda_{N+1}\rb/(\lambda_{N+1}-g_{i_0\mu_0\nu_0})
\end{array}$$
which satisfy $V_{i_0\mu_0\nu_0}\cap W_{i_0\mu_0\nu_0}=\emptyset$ and $V_{i_0\mu_0\nu_0}\cup W_{i_0\mu_0\nu_0}=X_{i_0\mu_0\nu_0}$. Hence it is enough to show that the lemma holds for both $V_{i_0\mu_0\nu_0}$ and $W_{i_0\mu_0\nu_0}$ instead of $X$. First we will find $PRD$s over $\Smn$ where the conclusion of the lemma holds for $V_{i_0\mu_0\nu_0}$. Define 
$$\begin{array}{rll}
G_{i_0\mu_0\nu_0}&=&\lambda^{\nu_0}+\underset{\nu\neq\nu_0,\mu}{\sum}a_{i_0\mu\nu}y_{M+1}y^{\mu}\lambda^{\nu}\\
&\equiv&y_{M+1}f_{i_0}\text{ mod }(y_{M+1}g_{i_0\mu_0\nu_0}-1)B_{i_0\mu_0\nu_0}\ls y_{M+1} \rs
\end{array}$$

Observe that  $G_{i_0\mu_0\nu_0}$ and $F=y_{M+1}g_{i_0\mu_0\nu_0}-1$ are preregular in $\lambda$ and $y$ respectively. Thus  after a Weierstrass change of variables among $\lambda$'s, $G_{i_0\mu_0\nu_0}$ becomes regular in $\lambda_N$, and after a change of variables among $y$'s, $F$ becomes regular in $y_{M+1}$. Note that both of these variable changes are automorphisms of $B_{i_0\mu_0\nu_0}\ls y_{M+1}\rs$. 
 
Applying the Weierstrass Division Theorem (Theorem 2.3.8 of \cite{RoSPS}) to first divide by $G_{i_0\mu_0\nu_0}$, and then by $F$, we see that for some ideal $I'$, there is a finite monomorphism
\begin{multline}\tag{$\ddag$}\label{psi}
\psi:A_{i_0\mu_0\nu_0}\ls y_1,...,y_M\rs\lb\lambda_1,...,\lambda_{N-1}\rb/I'\hookrightarrow B_{i_0\mu_0\nu_0}\ls y_{M+1} \rs/\\
(f_1,...,f_k, F).
\end{multline}

Lexicographically $(M,N-1)<(M,N)$, so applying the inductive hypothesis, there are finitely many parameterized families of R-domains $Y_j$ over $S_{m+M,n+N-1}$ defined by the parameter rings $A_j$ and the function rings $C_j$ (whose presentation do not involve $\lambda_N$ and $y_{M+1}$), and there are ideals $J_j\subset C_j $ satisfying $I'C_j\subset J_j$, integers $M_j\leq M$, $N_j\leq N-1$ and Weierstrass changes of variables  $\psi_j$ fixing the parameter rings $A_j$ such that each
$$ \psi_j:A_j/(J_j\cap A_j)\ls y_1,...,y_{M_j}\rs\lb\lambda_1,...,\lambda_{N_j}\rb\hookrightarrow C_j/J_j$$
is a finite monomorphism and also the assertions (i), (ii) and (iii) of the lemma are satisfied if we substitute $K'$-Dom$_{m+M,n+N-1}A_{i_0\mu_0\nu_0}\ls y_1,...,y_M\rs\lb\lambda_1,...,\lambda_{N-1}\rb$ for $\kDom B$,  $K'$-Dom$_{m+M,n+N-1}C_j$ for $\kDom B_j$, $I'$ for $I$ and $J_j$ for $I_j$. Furthermore by Lemma \ref{nil} we may assume that each $J_j$ is a radical ideal of $C_j$.

Now define 
$$B_j=C_j\lb \lambda_{N}\rb\ls y_{M+1}\rs/(y_{M+1}g_{i_0\mu_0\nu_0}-1)$$
and notice that each $B_j$ is a generalized ring of fractions over $\Smn$. Furthermore each $X_j=\text{Dom}_{\Smn}B_j$ is a $PRD$ over $\Smn$. 

Let us by $\bar{\psi}$ denote the change of variables that we do to obtain $\psi$ as in the Equation \eqref{psi} and let $\bar{\psi}(J_j)$ denote the set of elements we obtain through applying this operation on the ideal $J_j$. Let $I_j$ be the ideal of $B_j$ which is generated by $I\cup \bar{\psi}(J_j)$. With this notation Lemmas \ref{injective} and \ref{idealmove} imply that for each $j$ $\psi$ extends to a finite monomorphism
$$\bar{\psi}_j:C_j/J_j\rightarrow B_j/I_j$$
whose kernel is contained in the nilradical of $C_j/I_j$ which is trivial by the assumption that $J_j$ is a radical ideal of $C_j$. By Lemma \ref{flat} we also have $\bar{\psi}_j^{-1}(IB_j)\cap C_j\subset I'C_j\subset J_j$. Hence composing $\psi_j$ with $\bar{\psi}_j$ we get a finite monomorphism
$$\phi_j:A_j/(I_j\cap A_j)\ls y_1,...,y_{M_j}\rs\lb\lambda_1,...,\lambda_{N_j}\rb\hookrightarrow B_j/I_j,$$
as described in the statement (iv) of the lemma. In this case we also have $\kDom B_j\cap V(IB_j)_{K'}=\kDom B_j\cap V(I_j)_{K'}$ and the statements (i), (ii) and (iii) of the lemma follow easily with the substitution of $\kDom B_{i_0\mu_0\nu_0}\ls 1/g_{i_0\mu_0\nu_0}\rs$ for $\kDom B$.

For partitioning $W_{i_0\mu_0\nu_0}$ into $PRD$s, we follow the same lines as above, this time using
$$F=\lambda_{N+1}-g_{i_0\mu_0\nu_0}$$
which is preregular in $\lambda$. 

{\em Case 2.}  
\begin{multline*}B=A\ls y_1,...,y_{M+S}\rs\lb\lambda_1,...,\lambda_{N+T}\rb/(\{g_{M+i}y_{M+i}-h_{M+i}\}_{i=1}^S\cup\\
\{G_{N+l}\lambda_{N+l}-H_{N+l}\}_{l=1}^T).
\end{multline*}

Let $B'$ be $A\ls y_1,...,y_{M+S}\rs\lb\lambda_1,...,\lambda_{N+T}\rb$, then by the discussion following Lemma \ref{idealmove} there is an ideal $I'\subset B'$ corresponding to $I$. By Case 1 we can find finitely many $PRD$s $Y_j$ over $S_{m+M+S,n+N+T}$ defined by function rings $B'_j$, parameter rings $A_j$ such that there are ideals $I_j\subset B'_j$ containing $I'B_j$, integers $M_j$, $N_j$  and  finite monomorphisms $\phi_j$
$$\phi_j:A_j/(I'_j\cap A_j)\ls y_1,...,y_{M_j}\rs\lb\lambda_1,...,\lambda_{N_j}\rb\hookrightarrow B'_j/I'_j,$$
and the statements (i), (ii) and (iii) of the lemma are satisfied if we substitute $B'$ for $B$, $I'$ for $I$, $B'_j$ for $B_j$, $I'_j$ for $I_j$ and $K'$-Dom$_{m+M+S,n+N+T}$ for $\kDom$.

Notice that each $B'_j/I'_j$ can be thought of as a quotient of a generalized ring of fractions $B_j$ over $\Smn$ and an ideal $I_j$ of $B_j$. Notice also that the domain $X_j$ of each $B_j$ over $\Smn$ is a $PRD$. These $B_j$ and $I_j$ satisfy the statement (i), (ii) and (iii) of the lemma. 
\end{proof}

\begin{rem}\label{norm1r} 

i) This proof of Lemma \ref{norm1} still works if we allow $X$ to be the domain of an arbitrary generalized ring of fractions. In that case the domains $X_i$ are not necessarily $PRD$s, but they are domains of generalized rings of fractions.

ii) Although we did not obtain a normalization theorem for the algebras in consideration, the algebras $B_j/I_j$ that come up at the end have the nice property that the union of the associated $D$-semianalytic sets $\kDom B_j\cap V(I_j)_{K'}$ is the $D$-semianalytic set $\kdom B\cap V(I)_{K'}$ that we started with. Notice that the normalization process goes over the non-parameter $y$ and $\lambda$ variables as long as it is possible, so after the Weierstrass changes of variables $y_1,...,y_{M_j}$, $\lambda_1,...,\lambda_{N_j}$ become free variables. That is to say that intuitively the normalization process terminates when the ``relations'' remaining are only among the parameter variables $x$ and $\rho$.  
\end{rem}

In fact we can say something more about the Krull dimensions after the normalization as the following lemma indicates.

\begin{lem}\label{nicenorm}
Let $B$ be as in Lemma \ref{norm1} and suppose that $I$ is an ideal of $B$ satisfying $I=\mathcal{I}(\Dom B\cap V(I)_{K})$ then in Lemma \ref{norm1}, $I_j$ can be chosen such that k-dim $B/I\geq$ k-dim $B_j/I_j$. Moreover, in that case, for each minimal prime divisor $\mfp$ of $I_j$ there is a maximal ideal $\mfm\in\text{Dom}_{\Smn} B_j$ such that $\mfp\subset\mfm$. 

\end{lem}

\begin{proof}

 Apply Lemma \ref{norm1} to $I$, $B$, and $A$ to get $PRD$s $X_j$ defined by the rings of functions $B_j$ and parameter rings $A_j$, ideals $I_j\subset B_j$ and integers $M_j$, $N_j$ so that each
$$\phi_j:A_j/(I_j\cap A_j)\ls y_1,...,y_{M_j}\rs\lb\lambda_1,...,\lambda_{N_j}\rb\hookrightarrow B_j/I_j$$
is a finite monomorphism, each $\phi_j$ is a Weierstrass maps among $y$ and $\lambda$ separately and conditions (i), (ii) and (iii) of Lemma \ref{norm1} are satisfied. Let us fix a $j$ and write
$$\begin{array}{ccl}
A&=&S_{s,t}/(\{g_kx_k-h_k\}_{k=m+1}^s\cup\{\bar{g}_i\rho_i-\bar{h}_i\}_{i=n+1}^t)\\
A_j&=&S_{s',t'}/(\{g_kx_k-h_k\}_{k=m+1}^{s'}\cup\{\bar{g}_i\rho_i-\bar{h}_i\}_{i=n+1}^{t'})\\
B&=&(A\otimes_A^s S_{S,T})/(\{G_ky_k-H_k\}_{k=M+1}^S\cup\{\bar{G}_i\lambda_i-\bar{H}_i\}_{i=N+1}^T)\\
B_j&=&(A_j\otimes_A^sS_{S',T'})/(\{G_ky_k-H_k\}_{k=M+1}^{S'}\cup\{\bar{G}_i\lambda_i-\bar{H}_i\}_{i=N+1}^{T'})\\
\end{array}$$
where $m\leq s\leq s'$, $n\leq t\leq t'$ and $M\leq S\leq S'$, $N\leq T\leq T'$.

All is clear if $I_j=B_j$, hence assume that $I_j$ is a proper ideal of $B_j$ and let $J$ be  a radical ideal of $A_j$ which contains $I_j\cap A_j$. Notice that $\phi_j$ induces a finite map
$$\phi_j':A_j/J\ls y_1,...,y_{M_j}\rs\lb\lambda_1,...,\lambda_{N_j}\rb\hookrightarrow B_j/(I\cup J),$$
and by Lemma \ref{injective} the nilradical of $A_j/J\ls y_1,...,y_{M_j}\rs\lb\lambda_1,...,\lambda_{N_j}\rb$ contains the kernel of $\phi'_j$. By Lemma \ref{nil} that nilradical is the ideal $J$ and by Lemma \ref{idealmove} it follows that $\phi_j'$ is a finite monomorphism. Therefore by replacing $I_j$ with $I_j\cup \mathcal{I}(\dom A_j \cap V(I_j\cap A_j)_{K})\cdot B_j$ we may assume that $I_j\cap A_j=\mathcal{I}(\dom A_j \cap V(I_j\cap A_j)_{K})$.

Now assume that $\mfp$ is a minimal prime divisor of $I_j$ and let us write $\bar{\mfp}$ for the corresponding ideal in $S_{s'+S',t'+T'}$. Going back to the process of normalization we see that $H_k=1$ for $k>S$ and $\bar{G}_i=1$ for $i>T$, hence $G_k$ and $\bar{G}_i$ are not contained in $\bar{\mfp}$ for $k>S$ and $i>T$. On the other hand, by the Going Up Theorem (Theorem 9.4 of \cite{matsu}) for integral extensions, $\mfp\cap \phi_j(A_j/(I_j\cap A_j)\ls y_1,...,y_{M_j}\rs\lb\lambda_1,...,\lambda_{N_j}\rb)$ has to be a minimal prime of the ring $A_j/(I_j\cap A_j)\ls y_1,...,y_{M_j}\rs\lb\lambda_1,...,\lambda_{N_j}\rb$ and therefore by Lemma \ref{flat} $\bar{\mfp}\cap S_{s',t'}$ is a minimal prime divisor of the ideal corresponding to $I_j\cap A_j$ in $S_{s',t'}$. Because we have $I_j\cap A_j=\mathcal{I}(\dom A_j\cap V(I_j\cap A_j)_K)$, $\bar{\mfp}\cap S_{s',t'}$ can not contain $g_k$, $\bar{g}_i$ for $k>m$ and $i>n$, and therefore neither can $\bar{\mfp}$. Hence there is a maximal ideal $\mfm$ of $S_{s'+S',t'+T'}$ which contains $\bar{\mfp}$ but none of the $g_k$, $\bar{g}_i$ for $k>m$ and $i>n$. So by Lemma \ref{dimeq} we have 
$$\text{k-dim }B/I\geq\text{k-dim }(S_{s'+S',t'+T'})_{\mfm}/\bar{\mfp}_{\mfm}=\text{k-dim }B_j/\mfp.$$

Finally let $\mfp'$ be $\bar{\mfp}\cap S_{s+S,t+T}$. If there is a prime ideal $\bar{\mfq}$ between the ideal corresponding to $I$ in $S_{s+S,t+T}$ and $\mfp'$, and if $\mfq$ is the canonical image of $\bar{\mfq}$ in $B$, and $\mfq'$ is the ideal of $S_{s'+S',t'+T'}$ which corresponds to $\mfq B_j$ then by the argument above the $g_k$, $\bar{g}_i$  are not contained in $\mfq'$ for $k>m$ and $i>n$. Then again by Lemma \ref{dimeq} 
$$\text{k-dim }B/\mfq=\text{ k-dim }B_j/\mfq B_j<\text{k-dim }B_j/\mfp$$
contradicting the assumption that $\mfp$ is minimal over $I_j$ as $I\subset\mfq B_j$. Therefore $\mfp'$ is a minimal prime divisor of the ideal corresponding to $I$ in $S_{s+S,t+T}$ and hence it can not contain $G_k$, $\bar{G}_i$ for $M<k\leq S$ and $N<i\leq T$ by the assumption that $I=\mathcal{I}(\Dom B\cap V(I)_{K})$, and so neither can $\bar{\mfp}$. This proves the second statement of the lemma.
\end{proof}

\begin{rem}\label{nicenormr}
i) Notice that as in part (i) of Remark \ref{norm1r} the assumption that $\Dom B$ is a $PRD$ is again not necessary in the proof of the above lemma.

ii) In fact we can make another improvement in the statement of Lemma \ref{norm1}. We claim that each $I_j$ can be chosen such that 
$$I_j=\mathcal{I}(\Dom B_j\cap V(I_j)_{K}).$$ 

We can assume that $I=\mathcal{I}(\Dom B\cap V(I)_{K})$ and apply Lemma \ref{norm1}. Let the ideals $I_j$ be as in the statement of Lemma \ref{nicenorm} so that k-dim $B/I\geq$ k-dim $B_j/I_j$ for all $j$. It is clear that the statement is true if k-dim $B/I=0$. Assume that the claim holds for all generalized rings of fractions $C$ and ideals $J\subset C$ such that k-dim $C/J<d$ and assume k-dim $B/I=d$. Put
$$J_j=\mathcal{I}(\Dom B_j\cap V(I_j)_{K}),$$
and replace each $I_j$ with $J_j$ in the statement of Lemma \ref{norm1}. Observe that the assertions (i), (ii) and (iii) of Lemma \ref{norm1} still hold and the assertion (iv) is still true for all $J_j$ satisfying  k-dim $B_j/J_j=$ k-dim $B_j/I_j$. On the other hand for $J_j$ satisfying k-dim $B_j/J_j<$ k-dim $B_j/I_j\leq d$, we apply Lemma \ref{norm1} once again replacing $I$ with $J_j$ and $B$ with $B_j$. By induction and Lemma \ref{nicenorm} we have the claim.

\end{rem}

In the next lemma we take the normalization one step further by applying the normalization process also to the parameter rings $A_j$ of Lemma \ref{norm1}.

\begin{cor}\label{norm2}
Using the notation of Definition \ref{PFOR} let $X$ be a $PRD$ over $\Smn(E,K)$ defined by the function ring  $B$. Let  $I$ be  an ideal of $B$, then there exist finitely many $PRD$s $X_j$ defined by the function rings $B_j$, parameter rings $A_j$, together with ideals $I_j\subset B_j$ satisfying $IB_j\subset I_j$ and integers $m_j$, $n_j$, $M_j$, $N_j$, $m_j+n_j\leq m+n$, $M_j+N_j\leq M+N$ such that for any complete extension $K'$ of $K$

i) $\kDom B=\bigcup_{j}\kDom B_j$,

ii) $\kDom B\cap V(I)_{K'}=\bigcup_j (\kDom B_j\cap V(I_j)_{K'})$, 

iii) $\kDom B_j\cap V(I_j)_{K'}\cap \kDom B_i=\emptyset\text{ for }i\neq j,$ 

iv) writing each $B_j$ as 
\begin{multline*} B_j=A_j\ls y_1,....,y_{M+S_j}\rs\lb\lambda_1,...,\lambda_{N+T_j}\rb/(\{y_{M+k}g_{j,M+k}-h_{j,M+k}\}_{k=1}^{S_j}\cup\\
\{\lambda_{N+l} G_{j,M+l}-H_{j,M+l}\}_{l=1}^{T_j}),
\end{multline*}
where $A_j=S_{s_j,t_j}/J_j$ for some $J_j\in S_{s_j,t_j}/J_j$ and $S_{s_j,t_j}$ is the ring of separated power series over the variables $x_1,...,x_{s_j}$ and $\rho_1,...,\rho_{t_j}$ we have finite monomorphisms 
$$\phi_j:S_{m_j,n_j}\ls y_1,...,y_{M_j}\rs\lb\lambda_1,...,\lambda_{N_j}\rb\hookrightarrow B_j/I_j$$
such that 
$$\phi_j'=\phi_j|_{S_{m_j,n_j}}:S_{m_j,n_j}\hookrightarrow A_j/(I_j\cap A_j)$$
is also a finite monomorphism and $\phi_j$ is a Weierstrass change of variables among $x$, $\rho$, $y$ and $\lambda$ variables separately.
\end{cor}

\begin{proof}
After applying Lemma \ref{norm1} to $B$ and $I$, we apply it once again to $A_j$ and $I_j$ by considering $x$ and $\rho$ variables as non-parameter variables and apply Lemma \ref{injective}. The statement $m_j+n_j\leq m+n$, $M_j+N_j\leq M+N$ follows from Lemma \ref{nicenorm}.
\end{proof}

\section{Applications of Parameterized Normalization}

In this section we apply the normalization process of the previous section to prove our main results. We are primarily interested in establishing the relations between the dimensions that can be associated with $D$-semianalytic sets. See Definitions \ref{Adim}, \ref{gdim} and \ref{wdim} for these concepts of dimension. We establish equalities among these dimensions in Theorem \ref{alleq0}, Corollary \ref{locglob} and Lemma \ref{charpdim}. As a final application, we prove Theorem \ref{dimadd} which links the dimension of a $D$-semianalytic set over a field of characteristic $0$ with the dimension of its projection onto a coordinate hyperplane and dimensions of fibers of points in this projection.   

\begin{lem}\label{rdimdomwdim}
Let $B$ be a generalized ring of fractions over $\Smn$, $I$ an ideal of $B$ satisfying $I=\mathcal{I}(\Dom B\cap V(I)_K)$, and $\pi:(K^{\circ})^{m+M}\times(K^{\circ\circ})^{n+N}\rightarrow (K^{\circ})^m\times(K^{\circ\circ})^n$ the projection map. If $\pi(\Dom B\cap V(I)_K)$ is a somewhere dense set, then k-dim $B/I=$ dim$_{\Smn} B/I\geq m+n$.
\end{lem}
\begin{proof}
Apply Lemma \ref{nicenorm} to $B$ and $I$ considering $(K^{\circ})^m\times(K^{\circ\circ})^n$ to be the space of parameters to obtain generalized rings of fractions $B_j$, $A_j$, ideals $I_j\subset B_j$ and integers $M_j$, $N_j$ as described in the lemma. By Remark \ref{nicenormr} we may assume that $I_j=\mathcal{I}(\Dom B_j\cap V(I_j)_K)$. Notice that for some $j_0$, $\dom A_{j_0}\cap V(I_{j_0}\cap A_{j_0})_K$ has to be somewhere dense therefore by Lemma \ref{wdimgrdimeq1} k-dim $A_{j_0}/(I_{j_0}\cap A_{j_0})=$ dim$_{\smn}A_{j_0}/(I_{j_0}\cap A_{j_0})=m+n$, hence k-dim $B/I\geq$ k-dim $B_{j_0}/I_{j_0}\geq m+n$.
\end{proof}

Combining this lemma with Theorem \ref{stratify}, we have established the following.

\begin{thm}\label{alleq0}
Assume {\em Char }$K=0$. For a generalized ring of fractions $B$ over $\smn$ and ideal $I$ of $B$ satisfying $I=\mathcal{I}(\dom B\cap V(I)_K)$,
$$\text{w-dim }\dom B\cap V(I)_K=\text{g-dim }\dom B\cap V(I)_K=\text{ k-dim }B/I.$$
\end{thm}

The following result about local and global dimensions of a $D$-semianalytic set also easily follows.

\begin{cor}\label{locglob}
Assume {\em Char} $K=0$, then for a $D$-semianalytic set $X$ we have
$$\text{g-dim }X=\sup_{\pb\in X}\text{ g-dim }X_{\pb}.$$
\end{cor}
\begin{proof}
By the equality of the weak dimension and the geometric dimension in {\em Char }$K=0$, it is enough to consider sets of type $\dom B\cap V(I)_K$ where $B$ is a generalized ring of fractions over $\smn$ and $I=\mathcal{I}(\dom B\cap V(I)_K )$. By Theorem \ref{stratify}, $\dom B\cap V(I)$ contains a $(\text{dim}_{\smn}B/I)$-dimensional manifold and the result follows from Lemma \ref{dimeq}.
\end{proof}

 A weaker result in the case {\em Char }$K=p>0$ is the following.

\begin{lem}\label{charpdim}
Suppose that {\em Char }$K=p>0$, and $K$ is algebraically closed or $[K:K^p]<\infty$ and there is a countable basis for the topology of $K^{\circ}$. Then for a $D$-semianalytic set $X$ we have  
$$\text{g-dim } X= \sup_{\pb\in X}\text{ g-dim }X_{\pb}.$$ 
\end{lem}

\begin{proof}
 The case where $K$ is algebraically closed is proved in Lemma 2.3 of \cite{RoSPS}. In the case $[K:K^p]<\infty$ and there is a countable basis for the topology of $K^{\circ}$, by induction and Lemma \ref{dimeq} we see that for each $\pb\in X$ there is an $R$-domain $U$ such that $\pb$ is an element of the set $W$ of $K$-rational points in $U$ and there is an ideal $J\subset \mathcal{O}_K(U)$ such that for each $R$-domain $U'\subset U$ such that $\pb$ is an element of the set $W'$ of $K$-rational points in $U'$, k-dim $\mathcal{O}_K(U')/\mathcal{I}(W'\cap X)=$ k-dim $\mathcal{O}_K(U)/J$. Now by Theorem \ref{stratifyp} there is a neighborhood $W''$ of $\pb$ such that $W\cap W''\cap X$ is a union of $D$-semianalytic $K$-manifolds $X_i$ such that $\sup_i$ g-dim $X_i=$ k-dim $\mathcal{O}_K(U)/J$, which dominates w-dim $W\cap W''\cap X$ by Lemma \ref{rdimdomwdim}. Now the result follows from the assumption that there is a countable basis for the topology and the Baire Category Theorem for complete metric spaces.
\end{proof}

Next we will concentrate on the case {\em Char }$K=0$ and will prove additivity of dimensions of fibers and projections. We will need some groundwork before we can get to the main result.

For a subset $S$ of $(K^{\circ})^{m+M}\times (K^{\circ\circ})^{n+N}$ we will adopt the notation
\begin{multline*}
S^{(d)}:=\{\pb\in (K^{\circ})^{m}\times (K^{\circ\circ})^{n}: \text { w-dim }S(\pb)\geq d \text{ where } S(\pb) \text{ denotes}\\ 
\text{the fiber of the point }\pb\text{ in }S\}.
\end{multline*}

We are interested in the situation where $S^{(d)}$ is a somewhere dense subset of $(K^{\circ})^m\times (K^{\circ\circ})^n$. The next lemma shows that if $S$ has this property and is broken up into finitely many subsets, then one of those subsets also has this property.

\begin{lem}\label{propertyp}
 Let $S$ be a subset of $(K^{\circ})^{m+M}\times (K^{\circ\circ})^{n+N}$ and $S^{(d)}$ be somewhere dense. Suppose $S=\bigcup_{i=1}^kS_i$, then there is an $i_0$ such that $S_{i_0}^{(d)}$ is somewhere dense.
\end{lem}

\begin{proof}
Let $\pb\in S$, so that the fiber $S(\pb)$ of $\pb$ in $S$ is the union of fibers $S_{i}(\pb)$ of $\pb$ in each $S_i$. As there are only finitely many coordinate hyperplanes, we may assume that there is a $d$-dimensional coordinate projection $\pi$ such that $S^{(d)}$ consists of points $\pb$ such that the projection $\pi(S(\pb))$ are somewhere dense. Notice that $\pi(S(\pb))=\bigcup_{i=1}^k\pi(S_{i}(\pb))$ hence for each $\pb\in S^{(d)}$ there must be some $i_{\pb}$ such that $\pi(S_{i_{\pb}}(\pb))$ is somewhere dense. Therefore $\pb\in S^{(d)}$ implies $\pb\in S_{i_{\pb}}^{(d)}$ for some $i_{\pb}$, which in turn implies $S^{(d)}=\bigcup_{i=1}^k S_i^{(d)}$. So there must be an $i_0$ such that $S_{i_0}^{(d)}$ is somewhere dense.
\end{proof}

If $X\subset (K^{\circ})^{m+M}\times (K^{\circ\circ})^{n+N}$ is of geometric dimension $m+n+d$ then clearly for some choice of coordinate hyperplane $(K^{\circ})^{m}\times (K^{\circ\circ})^{n}$, $X^{(d)}$ contains a non-empty open set hence is somewhere dense. A natural question is the following. If we know $X^{(d)}$ to be somewhere dense, is it true that g-dim $X\geq m+n+d$? In characteristic $0$ the answer is positive.

\begin{thm}\label{dimadd}
With the above notation suppose {\em Char }$K=0$ and let $X$ be a $D$-semianalytic subset of $(K^{\circ})^{m+M}\times (K^{\circ\circ})^{n+N}$. If $X^{(d)}$ is somewhere dense, then g-dim $X\geq m+n+d$.
\end{thm}

\begin{proof}

Since a $D$-semianalytic set is a finite union of sets of the form 
$$\Dom B_i\cap V(I_i)_K,$$  
where each $B_i$ is a generalized ring of fractions over $\Smn$, by Lemma \ref{propertyp} we may assume that $X=\Dom B\cap V(I)_K$ for some generalized ring of fractions $B$ and an ideal $I=\mathcal{I}(\Dom B\cap V(I)_K)$. Notice that in this case the hypothesis for the second assertion of Lemma \ref{nicenorm} is satisfied and by Remark \ref{nicenormr}, Lemma \ref{nicenorm} and  Lemma \ref{propertyp} we can find generalized rings of fractions $B'$ over $\Smn$, $A'$ over $\smn$; an ideal $J\subset B'$ containing $IB'$, integers $M'$, $N'$ and a Weierstrass change of variables $\bar{\phi}$, such that $(\Dom B' \cap V(J)_K)^{(d)}$ is somewhere dense, $\Dom B\cap V(I)_K$ contains $\Dom B' \cap V(J)_K$, and $\bar{\phi}$ induces a finite monomorphism 
$$\phi:A'/(J\cap A')\ls y_1,...,y_{M'}\rs\lb\lambda_1,...,\lambda_{N'}\rb\hookrightarrow B'/J.$$

Furthermore we have 
$$\text{dim}_{\Smn}B/I=\text{ k-dim }B/I\geq\text{ k-dim }B'/J=\text{ dim}_{\Smn}B'/J.$$
By part (ii) of Remark \ref{nicenormr}, we can also assume that 
$$J=\mathcal{I}(\Dom B'\cap V(J)_K).$$

Next we observe that projection of $\Dom B'\cap V(J)_K$ onto $(K^{\circ})^m\times(K^{\circ\circ})^n$ is a somewhere dense set contained in $\dom A' \cap V(J\cap A')_K$. By Lemma \ref{wdimgrdimeq1}, we have k-dim $A'/(J\cap A')\geq$ dim$_{\smn}A'/(J\cap A')=m+n$. Observe also that  $\phi$ is an injection fixing parameters and taking fibers to fibers, therefore if $\pb\in (\Dom B'\cap V(J)_K)^{(d)}$ and $\mfm$ is the maximal ideal of $A'$ corresponding to $\pb$ then $(A'/(J\cap A')\ls y_1,...,y_{M'}\rs\lb \lambda_1,...,\lambda_{N'}\rb)/\mfm$, which is isomorphic to $S_{M',N'}$, is mapped injectively and finitely into $(B'/J)/\mfm$. Notice that $(B'/J)/\mfm$ is a quasi-affinoid algebra which defines the $D$-semianalytic subset $X(\pb)$ of $(K^{\circ})^M\times(K^{\circ\circ})^N$. By assumption w-dim $X(\pb)\geq d$ and hence $M'+N'\geq$ dim$_{S_{M,N}}(B'/J)/\mfm\geq d$ by Theorem \ref{alleq0}. Therefore
$$\text{k-dim }A'/(J\cap A')\ls y_1,...,y_{M'}\rs\lb\lambda_1,...,\lambda_{N'}\rb\geq m+n+d$$ 
and the statement follows from Theorem \ref{alleq0}.
\end{proof}

\bibliographystyle{asl}
\nocite{*}
\bibliography{dimtheory}

\end{document}